\documentclass[letterpaper,dvipdf]{article}
\usepackage{amssymb}
\usepackage{amsmath}
\usepackage{theorem}
\newtheorem{thm}{Theorem}[section]
\newtheorem{lem}{Lemma}[section]
\theorembodyfont{\rmfamily}
\newtheorem{exm}{Example}[section]
\makeatletter

\@addtoreset{equation}{section}
\makeatother
\usepackage[dvipdfmx]{graphicx}
\usepackage{dcolumn}
\def\al{{\alpha}}
\def\be{{\beta}}
\def\ga{{\gamma}}
\def\de{{\delta}}
\def\ep{{\varepsilon}}
\def\ka{{\kappa}}
\def\la{{\lambda}}
\def\si{{\sigma}}
\def\th{{\theta}}

\def\De{{\varDelta}}
\def\bxi{{\text{\boldmath $\xi$}}}
\def\bga{{\text{\boldmath $\ga$}}}

\def\Ga{{\varGamma}}

\def\bXi{{\text{\boldmath $\varXi$}}}

\def\u{{\text{\boldmath $u$}}}

\def\z{{\text{\boldmath $z$}}}

\def\Q{{\text{\boldmath $Q$}}}
\def\U{{\text{\boldmath $U$}}}

\def\X{{\text{\boldmath $X$}}}
\def\Y{{\text{\boldmath $Y$}}}
\def\Z{{\text{\boldmath $Z$}}}
\def\Nc{{\cal N}}
\def\infi{{\infty}}
\def\dd{\,{\rm d}}
\def\Rb{\mathbb{R}}
\def\non{{\nonumber}}
\newcommand{\qed}{\hfill\hbox{\rule{7pt}{7pt}}}
\begin{document}
\title{Estimation in a simple linear regression model with measurement error}
\author{Hisayuki Tsukuma\footnote{Faculty of Medicine, Toho University, 5-21-16 Omori-nishi, Ota-ku, Tokyo 143-8540, Japan, E-Mail: tsukuma@med.toho-u.ac.jp}}
\maketitle 

\begin{abstract}
This paper deals with the problem of estimating a slope parameter in a simple linear regression model, where independent variables have functional measurement errors.
Measurement errors in independent variables, as is well known, cause biasedness of the ordinary least squares estimator.
A general procedure for the bias reduction is presented in a finite sample situation, and some exact bias-reduced estimators are proposed.
Also, it is shown that certain truncation procedures improve the mean square errors of the ordinary least squares and the bias-reduced estimators.

\par\vspace{10pt}\noindent
{\it AMS 2010 subject classifications:} 
Primary 
62F10
; secondary 
62J07
.

\par\vspace{4pt}\noindent
{\it Key words and phrases:}
Bias correction,
errors-in-variables model,
functional relationship, 
mean square error, 
multivariate calibration problem, 
repeated measurement, 
shrinkage estimator, 
statistical control problem, 
statistical decision theory, 
structural relationship.
\end{abstract}
\section{Introduction}
\label{sec:intro}

Linear regression model with measurement errors in independent variables is of practical importance, and many theoretical and experimental approaches have been studied extensively for a long time.
Adcock (1877, 1878) first treated estimation of the slope in a simple linear measurement error model and derived the maximum likelihood (ML) estimator, which nowadays is known as orthogonal regression estimator (see Anderson (1984)).
Reiers\o{}l (1950) has investigated identifiability related to possibility of constructing a consistent estimator.
For efficient estimation, see Bickel and Ritov (1987) and, for consistent estimation based on shrinkage estimators, see Whittemore (1989) and Guo and Ghosh (2012).
A multivariate generalization of univariate linear measurement error model has been considered by Gleser (1981).
See Anderson (1984), Fuller (1987) and Cheng and Van Ness (1999) for a systematic overview of theoretical development in estimation of linear measurement error models.

Even though many estimation procedures for the slope have been developed and proposed, each procedure generally has both theoretical merits and demerits.
The ML estimator possesses consistency and asymptotic normality.
However, the first moment of the ML estimator does not exist and it is hard to theoretically investigate finite-sample properties of the ML procedure.
Besides the ML procedure, the most well-known procedure may be the least squares (LS) procedure.
The ordinary LS estimator has finite moments up to some order, but is not asymptotically unbiased.
The asymptotic biasedness of the LS estimator is called attenuation bias in the literature (see Fuller (1987)).

This paper addresses a simple linear measurement error model in a finite sample setup, and discusses the problem of reducing the bias and the mean square error (MSE) for slope estimators.
Suppose that the $Y_i$ and the $X_{ij}$ are observable variables for $i=1,\ldots,n$ and $j=1,\ldots,r$, where $r$ is the number of groups and $n$ is the sample size of each group.
Suppose also that the $Y_i$ and the $X_{ij}$ have the following model:
\begin{equation}\label{eqn:model0}
\begin{split}
Y_i&=\al_0+\be\ga_i+\de_i, \\
X_{ij}&=\ga_i+\ep_{ij},
\end{split}
\end{equation}
where $\al_0$ and $\be$ are, respectively, unknown intercept and slope parameters, the $\ga_i$ are unobservable latent variables, and the $\de_i$ and the $\ep_{ij}$ are random error terms.
Assume that the $\de_i$ and the $\ep_{ij}$ are mutually independent and distributed as $\de_i\sim\Nc(0,\tau^2)$ and $\ep_{ij}\sim\Nc(0,\si_x^2)$, respectively, where $\tau^2$ and $\si_x^2$ are unknown.
It is important to note that the error variance in independent variables, $\si_x^2$, can be estimated.

For the latent variables $\ga_i$ in model \eqref{eqn:model0}, there are two different points of view, namely, the $\ga_i$ are considered as unknown fixed values or as random variables.
In the former case, \eqref{eqn:model0} is referred to as a functional model and, in the latter case, is called a structural model (Kendall and Stuart (1979), Anderson (1984) and Cheng and Van Ness (1999)).
In this paper, we assume the functional model and shall develop a finite-sample theory of estimating the slope $\be$.

The remainder of this paper is organized as follows.
In Section \ref{sec:canonical}, we simplify the estimation problem in model \eqref{eqn:model0}, and define a broad class of slope estimators including the LS estimator, the method of moments estimator, and a Stefanski's (1985) estimator.
Also, Section \ref{sec:canonical} shows some technical lemmas used for evaluating moments.
Section \ref{sec:bias} presents a unified method of reducing the bias of the broad class as well as that of the LS estimator.
In Section \ref{sec:MSE}, we handle the problem of reducing the MSEs of slope estimators.
It is revealed that the slope estimation under the MSE criterion is closely related to the statistical control problem (see Zellner (1971) and Aoki (1989)) and also to the multivariate calibration problem (see Osborne (1991), Brown (1993) and Sundberg (1999)).
Our approach to the MSE reduction is carried out in a similar way to Kubokawa and Robert (1994), and a general method is established for improvement of several estimators such as the LS estimator and Guo and Ghosh's (2012) estimator.
Section \ref{sec:num} illustrates numerical performance for the biases and the MSEs of alternative estimators.
In Section \ref{sec:remarks}, we point out some remarks on our results and related topics.

\section{Simplification of the estimation problem}\label{sec:canonical}
\subsection{Reparametrized model}\label{subsec:reparametrize}

Define $\overline{X}_i=(1/r)\sum_{j=1}^r X_{ij}$ for $i=1,2,\ldots,n$.
Consider the regression of the $Y_i$ on the $\overline{X}_i$.
The LS estimator of $(\be,\al_0)$ is defined as a unique solution of
$$
\min_{\substack{-\infi<\be<\infi\\[1pt]-\infi<\al_0<\infi}}\ \sum_{i=1}^{n}(Y_i-\al_0-\be\overline{X}_i)^2.
$$
Denote by $(\hat{\be}^{LS},\hat{\al}_0^{LS})$ the resulting ordinary LS estimator of $(\be,\al_0)$.
Then $\hat{\be}^{LS}$ and $\hat{\al}_0^{LS}$ are given, respectively, by
$$
\hat{\be}^{LS}=\frac{\sum_{i=1}^{n}(\overline{X}_i-\overline{X})(Y_i-\overline{Y})}{\sum_{i=1}^{n}(\overline{X}_i-\overline{X})^2},\qquad
\hat{\al}_0^{LS}=\overline{Y}-\hat{\be}^{LS}\overline{X},
$$
where $\overline{X}=(1/n)\sum_{i=1}^{n}\overline{X}_i$ and $\overline{Y}=(1/n)\sum_{i=1}^{n}Y_i$.

Let $\bga=(\ga_1,\ldots,\ga_{n})^t$, $\Y=(Y_1,\ldots,Y_{n})^t$ and $\X=(\overline{X}_1,\ldots,\overline{X}_{n})^t$.
Define
$$
S=\frac{1}{r}\sum_{i=1}^{n}\sum_{j=1}^r(X_{ij}-\overline{X}_i)^2.
$$
Denote by $I_{n}$ the identity matrix of order $n$ and by $1_{n}$ the $n$-dimensional vector consisting of ones.
It is then observed that
\begin{equation}\label{eqn:model1}
\begin{aligned}
\Y&\sim\Nc_{n}(\al_01_{n}+\be\bga,\tau^2I_{n}), \\
\X&\sim\Nc_{n}(\bga,\si^2I_{n}),
\end{aligned}
\qquad
\begin{aligned}
&\\
S&\sim \si^2\chi_m^2,
\end{aligned}
\end{equation}
for $m=n(r-1)$ and $\si^2=\si_x^2/r$.
Note that $\Y$, $\X$ and $S$ are mutually independent.

Furthermore, let $\Q$ be an $n\times n$ orthogonal matrix whose first row is $1_{n}^t/\sqrt{n}$.
Denote $p=n-1$ and $\al=\al_0\sqrt{n}$.
Define $\Q\Y=(Z_0,\Z^t)^t$, $\Q\X=(U_0,\U^t)^t$ and $\Q\bga=(\th,\bxi^t)^t$, where $\Z$, $\U$ and $\bxi$ are $p$-dimensional vectors.
Then model \eqref{eqn:model1} can be replaced with
\begin{equation}\label{eqn:model2}
\begin{aligned}
Z_0&\sim\Nc(\al+\be\th,\tau^2),\\
U_0&\sim\Nc(\th,\si^2),
\end{aligned}
\qquad
\begin{aligned}
\Z&\sim\Nc_p(\be\bxi,\tau^2I_p),\\
\U&\sim\Nc_p(\bxi,\si^2I_p),
\end{aligned}
\qquad
\begin{aligned}
&\\
S&\sim \si^2\chi_m^2.
\end{aligned}
\end{equation}
These five statistics, $Z_0, \Z, U_0, \U$ and $S$, are mutually independent, and $\al$, $\be$, $\th$, $\bxi$, $\si^2$ and $\tau^2$ are unknown parameters.
Throughout this paper, we suppose that $\bxi\ne 0_p$.

From reparametrized model \eqref{eqn:model2}, the ordinary LS estimators $\hat{\be}^{LS}$ and $\hat{\al}^{LS}=\hat{\al}_0^{LS}\sqrt{n}$ can be rewritten, respectively, as
\begin{equation}\label{eqn:LS}
\hat{\be}^{LS}=\frac{\U^t\Z}{\Vert\U\Vert^2},\qquad
\hat{\al}^{LS}=Z_0-\hat{\be}^{LS}U_0.
\end{equation}

Hereafter, we mainly deal with the problem of estimating $\be$ in reparametrized model \eqref{eqn:model2}.
Denote the bias and the MSE of an estimator $\hat{\be}$, respectively, by
\begin{align*}
{\rm Bias}(\hat{\be};\be)&=E[\hat{\be}]-\be, \\
{\rm MSE}(\hat{\be};\be)&=E[(\hat{\be}-\be)^2],
\end{align*}
where the expectation $E$ is taken with respect to \eqref{eqn:model2}.
The bias of $\hat\be$ is smaller than that of another estimator $\hat\be_*$ if $|{\rm Bias}(\hat{\be};\be)|\leq|{\rm Bias}(\hat{\be}_*;\be)|$ for any $\be$.
Similarly, if ${\rm MSE}(\hat{\be};\be)\leq{\rm MSE}(\hat{\be}_*;\be)$ for any $\be$, then the MSE of $\hat\be$ is said to be better than that of $\hat\be_*$, or $\hat\be$ is said to dominate $\hat\be_*$.

\subsection{A class of estimators}

If $\lim_{n\to \infi} \Vert\bxi\Vert^2/p=\si_\xi^2$ where $\si_\xi^2$ is a positive value, it follows that $\U^t\Z/p\to \be\si_\xi^2$ and $\Vert\U\Vert^2/p\to \si_\xi^2+\si^2$ in probability as $n$ tends to infinity, and hence
\begin{equation}\label{eqn:attenuation}
\hat{\be}^{LS}\mathop{\to} \frac{\si_\xi^2}{\si_\xi^2+\si^2}\be \quad\text{in probability} \quad (n\to\infi).
\end{equation}
This implies that the ordinary LS estimator $\hat{\be}^{LS}$ is inconsistent and, more precisely, it is asymptotically biased toward zero.
This phenomenon is called attenuation bias (see Fuller (1987)).

For reducing the influence of attenuation bias, various alternatives to $\hat{\be}^{LS}$ have been proposed in the literature.
For example, a typical alternative is the method of moments estimator
\begin{equation}\label{eqn:MM}
\hat{\be}^{MM}=\frac{\U^t\Z/p}{\Vert\U\Vert^2/p-S/m}.
\end{equation}
The method of moments estimator $\hat{\be}^{MM}$ converges to $\be$ in probability as $n$ goes to infinity, but $\hat{\be}^{MM}$ does not have finite moments.
Noting that $\hat{\be}^{MM} =\{1-(p/m)S/\Vert\U\Vert^2\}^{-1}\hat{\be}^{LS}$ and also using the Maclaurin expansion $(1-x)^{-1}=\sum_{j=0}^\infi x^j$, we obtain the $\ell$-th order corrected estimator of the form
\begin{equation}\label{eqn:ST}
\hat{\be}_\ell^{ST}
=\bigg\{1+\frac{p}{m}\frac{S}{\Vert\U\Vert^2}+\cdots+\bigg(\frac{p}{m}\frac{S}{\Vert\U\Vert^2}\bigg)^\ell\, \bigg\}\hat{\be}^{LS}.
\end{equation}
The above estimator can also be derived from using the same arguments as in Stefanski (1985), who approached to the bias correction from Huber's (1981) M estimation.
However, it is still not known whether or not the bias of $\hat{\be}_\ell^{ST}$ is smaller than that of $\hat{\be}^{LS}$ in a finite sample situation.

Convergence \eqref{eqn:attenuation} is equivalent that $\bar{\be}=(1+\si^2/\si_\xi^2)\hat{\be}^{LS}$ converges to $\be$ in probability as $n$ goes to infinity.
Replacing $\si^2/\si_\xi^2$ of $\bar{\be}$ with a suitable function $\phi$ of $\Vert\U\Vert^2/S$ yields a general class of estimators,
\begin{equation}\label{eqn:class}
\hat{\be}_\phi=\bigg\{1+\phi\bigg(\frac{\Vert\U\Vert^2}{S}\bigg)\bigg\}\hat{\be}^{LS}.
\end{equation}
Note that $\hat\be^{MM}$ and $\hat\be_\ell^{ST}$ belong to the class \eqref{eqn:class}.
In this paper, we search a bias-reduced or an MSE-reduced estimator within \eqref{eqn:class} as an alternative to $\hat\be^{LS}$.

\subsection{Some useful lemmas}

Next, we provide some technical lemmas which form the basis for evaluating the bias and MSE of \eqref{eqn:class}.

\begin{lem}\label{lem:expectations}
Let $\U\sim\Nc_p(\bxi,\si^2I_p)$ and $S\sim \si^2\chi_m^2$.
Let $\phi$ be a function on the positive real line.
Define $\la=\Vert\bxi\Vert^2/(2\si^2)$ and denote by $P_\la(k)=e^{-\la}\la^k/k!$ the Poisson probabilities for $k=0,1,2,\ldots.$
Let $g_n(t)$ be the p.d.f.\ of $\chi^2_n$.
\begin{enumerate}
\item[{\rm(i)}]
If $E[|\phi(\Vert\U\Vert^2/S)\U^t\bxi|/\Vert\U\Vert^2]<\infi$ then we have
$$
E\bigg[\phi\Big(\frac{\Vert\U\Vert^2}{S}\Big)\frac{\U^t\bxi}{\Vert\U\Vert^2}\bigg]
=\sum_{k=0}^\infi \frac{2\la}{p+2k}P_\la(k)I_1(k|\phi),
$$
where $I_1(k|\phi)=\int_0^\infi\!\!\int_0^\infi\phi(w/s)g_{p+2k}(w)\dd w \,g_m(s)\dd s$.
\item[{\rm(ii)}]
If $E[|\phi(\Vert\U\Vert^2/S)|(\U^t\bxi)^2/\Vert\U\Vert^4]<\infi$ then we have
$$
E\bigg[\phi\Big(\frac{\Vert\U\Vert^2}{S}\Big)\frac{(\U^t\bxi)^2}{\Vert\U\Vert^{4}}\bigg]
=\sum_{k=0}^\infi \frac{2\la(1+2k)}{p+2k}P_\la(k)I_2(k|\phi),
$$
where $I_2(k|\phi)=\int_0^\infi\!\!\int_0^\infi w^{-1}\phi(w/s)g_{p+2k}(w)\dd w \,g_m(s)\dd s$.
\end{enumerate}
\end{lem}

When $\phi\equiv 1$, (i) and (ii) of Lemma \ref{lem:expectations} are, respectively,
\begin{align}
E\bigg[\frac{\U^t\bxi}{\Vert\U\Vert^2}\bigg]&=E\bigg[\frac{2\la}{p+2K}\bigg]\quad \textup{for $p\geq 2$},
\label{eqn:LS-1st}\\
E\bigg[\frac{(\U^t\bxi)^2}{\Vert\U\Vert^{4}}\bigg]
&=E\bigg[\frac{2\la(1+2K)}{(p+2K)(p+2K-2)}\bigg]\quad \textup{for $p\geq 3$},
\label{eqn:LS-2nd}
\end{align}
where $K$ is the Poisson random variable with mean $\la=\Vert\bxi\Vert^2/(2\si^2)$.
Identities \eqref{eqn:LS-1st} and \eqref{eqn:LS-2nd} have been given, for example, in Nishii and Krishnaiah (1988, Lemma 3).

\medskip
{\bf Proof of Lemma \ref{lem:expectations}.}\ \ 
(i) 
Denote
$$
E_1=E\bigg[\phi\Big(\frac{\Vert\U\Vert^2}{S}\Big)\frac{\U^t\bxi}{\Vert\U\Vert^2}\bigg].
$$
Let $\bxi_1=\bxi/\si$.
It turns out that
\begin{equation*}\label{eqn:ex1-1}
E_1=(2\pi)^{-p/2}\int_0^\infi\!\!\!\int_{\Rb^p}\phi\Big(\frac{\Vert\u\Vert^2}{s}\Big)\frac{\u^t\bxi_1}{\Vert\u\Vert^2}e^{-\Vert\u-\bxi_1\Vert^2/2}\dd\u\,g_m(s)\dd s.
\end{equation*}
Denote $c_0=(2\pi)^{-p/2}e^{-\la}$.
Let $\bXi$ be a $p\times p$ orthogonal matrix whose first row is $\bxi_1/\Vert\bxi_1\Vert$.
Making the orthogonal transformation $\u=(u_1,u_2,\ldots,u_p)^t\to\bXi^t\u$ gives that
\begin{equation}\label{eqn:ex1-2}
E_1
=c_0\int_0^\infi\!\!\!\int_{\Rb^p}\phi\Big(\frac{\Vert\u\Vert^2}{s}\Big)\frac{u_1\Vert\bxi_1\Vert}{\Vert\u\Vert^2}e^{-\Vert\u\Vert^2/2+u_1\Vert\bxi_1\Vert}\dd\u\,g_m(s)\dd s.
\end{equation}

Now, for $p\geq 2$, we make the following polar coordinate transformation
$$
\u
=\begin{pmatrix} u_1 \\ u_2 \\ u_3 \\ \vdots \\ u_{p-1} \\ u_p \end{pmatrix}
= \rho \left(\begin{array}{l}
\cos\varphi \\
\sin\varphi \cos\varphi_2 \\
\sin\varphi \sin\varphi_2 \cos\varphi_3 \\
\vdots \\
\sin\varphi \sin\varphi_2 \sin\varphi_3 \cdots \sin\varphi_{p-2}\cos\varphi_{p-1} \\
\sin\varphi \sin\varphi_2 \sin\varphi_3 \cdots \sin\varphi_{p-2}\sin\varphi_{p-1} \end{array}\right),
$$
where $\rho>0$, $0<\varphi<\pi$, $0<\varphi_i<\pi$ $(i=2,3,\ldots,p-2)$ and $0<\varphi_{p-1}<2\pi$.
The Jacobian of transformation $\u\to(\rho,\varphi,\varphi_2,\varphi_3,\ldots,\varphi_{p-1})$ is given by $\rho^{p-1}\sin^{p-2}\varphi \sin^{p-3}\varphi_2 \cdots \sin\varphi_{p-2}$, so \eqref{eqn:ex1-2} can be rewritten as
\begin{align*}\label{eqn:ex1-3}
E_1=c_1\int_0^\infi\!\!\!\int_0^\infi\!\!\!\int_0^\pi & \phi\Big(\frac{\rho^2}{s}\Big)
\frac{\Vert\bxi_1\Vert\cos\varphi}{\rho}e^{-\rho^2/2+\rho\Vert\bxi_1\Vert\cos\varphi}\non\\
& \times \rho^{p-1}\sin^{p-2}\varphi\dd\varphi\dd\rho\,g_m(s)\dd s,
\end{align*}
with
$$
c_1
=c_0\int_0^{2\pi}\dd\varphi_{p-1}\prod_{i=2}^{p-2}\int_0^\pi\sin^{p-i-1}\varphi_i\dd\varphi_i.
$$

Note here that, for an even $n$,
$$
\int_0^\pi \sin^m\varphi \cos^n\varphi \dd\varphi = 
\frac{\Ga[(m+1)/2]\Ga[(n+1)/2]}{\Ga[(m+n+2)/2]}
$$
and, for an odd $n$, the above definite integral is zero.
Thus, it is seen that
$$
c_1
=\frac{2^{1-p/2}\pi^{-1/2}e^{-\la}}{\Ga[(p-1)/2]}
$$
and
\begin{align*}
\int_0^\pi e^{\rho\Vert\bxi_1\Vert\cos\varphi} \cos\varphi \sin^{p-2}\varphi\dd\varphi
&=\sum_{j=0}^\infi \frac{\rho^j\Vert\bxi_1\Vert^j}{j!} \int_0^\pi \cos^{j+1}\varphi \sin^{p-2}\varphi\dd\varphi\\
&=\sum_{k=0}^\infi\frac{\rho^{2k+1}}{k!}\la^k\frac{\pi^{1/2}\Vert\bxi_1\Vert\Ga[(p-1)/2]}{2^k(p+2k)\Ga[(p+2k)/2]},
\end{align*}
so that
$$
E_1=\sum_{k=0}^\infi \frac{2\la}{p+2k}P_\la(k) \int_0^\infi\!\!\!\int_0^\infi \phi\Big(\frac{\rho^2}{s}\Big) \frac{\rho^{p+2k-1}e^{-\rho^2/2}}{\Ga[(p+2k)/2]2^{p/2+k-1}}\dd\rho\,g_m(s)\dd s.
$$
The change of variables $w=\rho^2$ leads to completeness of the proof of (i).

(ii)
Denote
$$
E_2=E\bigg[\phi\Big(\frac{\Vert\U\Vert^2}{S}\Big)\frac{(\U^t\bxi)^2}{\Vert\U\Vert^{4}}\bigg].
$$
Using the same arguments as in the proof of (i), we obtain
\begin{align*}
E_2=c_1\int_0^\infi\!\!\!\int_0^\infi\!\!\!\int_0^\pi& \phi\Big(\frac{\rho^2}{s}\Big)
\frac{\Vert\bxi_1\Vert^2\cos^2\varphi}{\rho^2}e^{-\rho^2/2+\rho\Vert\bxi_1\Vert\cos\varphi} \\
&\times \rho^{p-1}\sin^{p-2}\varphi\dd\varphi\dd\rho\,g_m(s)\dd s,
\end{align*}
Since
$$
\int_0^\pi e^{\rho\Vert\bxi_1\Vert\cos\varphi}\cos^2\varphi \sin^{p-2}\varphi\dd\varphi=\sum_{k=0}^\infi\frac{\rho^{2k}}{k!}\la^k\frac{\pi^{1/2}(1+2k)\Ga[(p-1)/2]}{2^k(p+2k)\Ga[(p+2k)/2]},
$$
it is observed that
$$
E_2=\sum_{k=0}^\infi P_\la(k) \frac{2\la(1+2k)}{p+2k} I_2(k|\phi),
$$
where
\begin{align*}
I_2(k|\phi)
&=\int_0^\infi\!\!\!\int_0^\infi\frac{1}{\rho^2}\phi\Big(\frac{\rho^2}{s}\Big)\frac{\rho^{p+2k-1}e^{-\rho^2/2}}{\Ga[(p+2k)/2]2^{p/2+k-1}}\dd\rho\,g_m(s)\dd s \\
&=\int_0^\infi\!\!\!\int_0^\infi \frac{1}{w}\phi\Big(\frac{w}{s}\Big)g_{p+2k}(w)\dd w \,g_m(s)\dd s.
\end{align*}
Hence the proof of (ii) is complete.
\qed

\begin{lem}\label{lem:chi}
Let $\U\sim\Nc_p(\bxi,\si^2I_p)$.
Let $i$ be a natural number such that $i<p/2$.
Denote by $K$ the Poisson random variable with mean $\la=\Vert\bxi\Vert^2/(2\si^2)$.
Then we have
$$
E\bigg[\frac{\si^{2i}}{\Vert\U\Vert^{2i}}\bigg]
=\begin{cases}
\prod_{j=1}^i(p-2j)^{-1} & \textup{if $\bxi=0_p$}, \\
E\big[\prod_{j=1}^i(p+2K-2j)^{-1}\big] & \textup{otherwise}.
\end{cases}
$$
\end{lem}

{\bf Proof.}\ \ 
We employ the same notation as in Lemma \ref{lem:expectations}.
Note that, when $\bxi\ne 0_p$, $\Vert\U\Vert^2/\si^2$ follows the noncentral chi-square distribution with $p$ degrees of freedom and noncentrality parameter $\Vert\bxi\Vert^2/\si^2$.
Since the p.d.f. of the noncentral chi-square distribution is given by $\sum_{k=0}^\infi P_\la(k)g_{p+2k}(w)$, it is seen that
\begin{align*}
E\bigg[\frac{\si^{2i}}{\Vert\U\Vert^{2i}}\bigg]
&=\sum_{k=0}^\infi P_\la(k)\int_0^\infi w^{-i}g_{p+2k}(w)\dd w \\
&=\sum_{k=0}^\infi P_\la(k)\prod_{j=1}^i(p+2k-2j)^{-1}
 =E\bigg[\prod_{j=1}^i(p+2K-2j)^{-1}\bigg] 
\end{align*}
for $p-2i>0$.
If $\bxi=0_p$, then $\Vert\U\Vert^2/\si^2\sim\chi_p^2$, so that $E[\si^{2i}/\Vert\U\Vert^{2i}]=\prod_{j=1}^i(p-2j)^{-1}$ for $p-2i>0$.
Thus the proof is complete.
\qed

\medskip
The following lemma is given in Hudson (1978).
\begin{lem}\label{lem:poisson}
Let $K$ be a Poisson random variable with mean $\la$.
Let $g$ be a function satisfying $|g(-1)| < \infty$ and $E[|g(K)|]<\infi$.
Then we have $\la E [g(K)]=E[Kg(K-1)]$.
\end{lem}

\section{Bias reduction}\label{sec:bias}

In this section, some results are presented for the bias reduction in slope estimation.
First, we give an alternative expression for the bias of the LS estimator $\hat{\be}^{LS}$.
\begin{lem}\label{lem:bias_LS}
Let $K$ be a Poisson random variable with mean $\la=\Vert\bxi\Vert^2/(2\si^2)$.
If $p\geq 2$, then the bias of $\hat{\be}^{LS}$ is finite.
Furthermore, if $p\geq 3$, the bias of $\hat{\be}^{LS}$ can be expressed as
\begin{equation*}
{\rm Bias}(\hat{\be}^{LS};\be)=-E\bigg[\frac{p-2}{p+2K-2}\bigg]\be.
\end{equation*}
\end{lem}

{\bf Proof.}\ \ 
Using identity \eqref{eqn:LS-1st} gives that for $p\geq 2$
\begin{equation}\label{eqn:bias-LS1}
{\rm Bias}(\hat{\be}^{LS};\be)=E\bigg[\frac{\U^t\bxi}{\Vert\U\Vert^2}\bigg]\be-\be
=E\bigg[\frac{2\la}{p+2K}\bigg]\be-\be.
\end{equation}
If $p\geq 3$, we apply Lemma \ref{lem:poisson} to \eqref{eqn:bias-LS1} so as to obtain
$$
{\rm Bias}(\hat{\be}^{LS};\be)
=E\bigg[\frac{2K}{p+2K-2}\bigg]\be-\be=-E\bigg[\frac{p-2}{p+2K-2}\bigg]\be.
$$
Hence the proof is complete.
\qed

\medskip
Let $\ell$ be a nonnegative integer.
Define a simple modification of $\hat{\be}_\ell^{ST}$, given in \eqref{eqn:ST}, as
\begin{equation}\label{eqn:BRi}
\hat{\be}_\ell^{BR}=\bigg\{1+\sum_{j=1}^\ell\frac{a_j}{b_j}\bigg(\frac{S}{\Vert\U\Vert^2}\bigg)^j\bigg\}\hat{\be}^{LS},
\end{equation}
where $a_j=(p-2)(p-4)\cdots(p-2j)$ and $b_j=m(m+2)\cdots(m+2j-2)$ for $j=1,\ldots,\ell$, and $\hat{\be}_0^{BR}\equiv \hat{\be}^{LS}$.
We then obtain the following lemma.
\begin{lem}\label{lem:bias-BRi}
Let $K$ be a Poisson random variable with mean $\la=\Vert\bxi\Vert^2/(2\si^2)$.
Assume that $p\geq 5$.
If $\ell<(p-2)/2$, then ${\rm Bias}(\hat{\be}_\ell^{BR};\be)$ can be expressed as
$$
{\rm Bias}(\hat{\be}_\ell^{BR};\be)=-E\bigg[\prod_{j=1}^{\ell+1}\frac{p-2j}{p+2K-2j}\bigg]\be.
$$
\end{lem}

{\bf Proof.}\ \ 
We prove a case when $\ell\geq 1$ because the $\ell=0$ case is equivalent to Lemma \ref{lem:bias_LS}.
Note that
\begin{equation*}
E[\hat{\be}_\ell^{BR}]=E[\hat{\be}^{LS}]+E\bigg[\sum_{j=1}^\ell\frac{a_j}{b_j}\bigg(\frac{S}{\Vert\U\Vert^2}\bigg)^j\frac{\U^t\bxi}{\Vert\U\Vert^2}\bigg]\be,
\end{equation*}
which implies from Lemma \ref{lem:bias_LS} that
\begin{equation}\label{eqn:EBR2}
{\rm Bias}(\hat{\be}_\ell^{BR};\be)=-E\bigg[\frac{p-2}{p+2K-2}\bigg]\be+E\bigg[\sum_{j=1}^\ell\frac{a_j}{b_j}\bigg(\frac{S}{\Vert\U\Vert^2}\bigg)^j\frac{\U^t\bxi}{\Vert\U\Vert^2}\bigg]\be.
\end{equation}

Since $E[X^j]=b_j$ for $j=1,\ldots,\ell$ when $X\sim\chi^2_m$, using (i) of Lemma \ref{lem:expectations} and Lemma \ref{lem:chi} gives 
\begin{align}\label{eqn:EBR3}
&E\bigg[\frac{a_j}{b_j}\bigg(\frac{S}{\Vert\U\Vert^2}\bigg)^j\frac{\U^t\bxi}{\Vert\U\Vert^2}\bigg]\non\\
&=\frac{a_j}{b_j}\sum_{k=0}^\infi \frac{2\la}{p+2k}P_\la(k)\int_0^\infi\!\!\!\int_0^\infi\frac{s^j}{w^j}g_{p+2k}(w)\dd w \,g_m(s)\dd s \non \\
&=a_j\sum_{k=0}^\infi \frac{2\la}{p+2k}P_\la(k)\int_0^\infi\frac{1}{w^j}g_{p+2k}(w)\dd w \non \\
&=E\bigg[\frac{2\la}{p+2K} \prod_{i=1}^{j}\frac{p-2i}{p+2K-2i}\bigg] 
\end{align}
for $p-2j>0$.
Applying Lemma \ref{lem:poisson} to \eqref{eqn:EBR3} gives that for $p-2-2j>0$
$$
E\bigg[\frac{a_j}{b_j}\bigg(\frac{S}{\Vert\U\Vert^2}\bigg)^j\frac{\U^t\bxi}{\Vert\U\Vert^2}\bigg]
=E\bigg[\frac{2K}{p+2K-2} \prod_{i=1}^{j}\frac{p-2i}{p+2K-2i-2}\bigg],
$$
which is substituted into \eqref{eqn:EBR2} to obtain
$$
{\rm Bias}(\hat{\be}_\ell^{BR};\be)=-E\bigg[\frac{p-2}{p+2K-2}-\frac{2K}{p+2K-2}\sum_{j=1}^\ell \prod_{i=1}^{j}\frac{p-2i}{p+2K-2i-2}\bigg]\be.
$$

It is here observed that
\begin{align*}
&\frac{p-2}{p+2K-2}-\frac{2K}{p+2K-2}\sum_{j=1}^\ell \prod_{i=1}^{j}\frac{p-2i}{p+2K-2i-2}\\
&=\prod_{j=1}^{2}\frac{p-2j}{p+2K-2j}-\frac{2K}{p+2K-2}\sum_{j=2}^\ell \prod_{i=1}^{j}\frac{p-2i}{p+2K-2i-2}\\
&=\cdots=\prod_{j=1}^{\ell+1}\frac{p-2j}{p+2K-2j},
\end{align*}
which yields that, for $p-2\ell-2>0$, 
$$
{\rm Bias}(\hat{\be}_\ell^{BR};\be)=-E\bigg[\prod_{j=1}^{\ell+1}\frac{p-2j}{p+2K-2j}\bigg]\be.
$$
Hence the proof is complete.
\qed

\begin{exm}
If $k$ is a nonnegative integer and $\ell\geq 1$, it follows that
$$
0<
\prod_{j=1}^{\ell+1}\frac{p-2j}{p+2k-2j}
\leq \prod_{j=1}^\ell\frac{p-2j}{p+2k-2j}
\leq\cdots\leq \frac{p-2}{p+2k-2}.
$$
Combining Lemmas \ref{lem:bias_LS} and \ref{lem:bias-BRi} immediately yields that, for any $\be$,
$$
|{\rm Bias}(\hat{\be}_\ell^{BR};\be)|
\leq |{\rm Bias}(\hat{\be}_{\ell-1}^{BR};\be)|
\leq\cdots
\leq |{\rm Bias}(\hat{\be}_1^{BR};\be)|
\leq |{\rm Bias}(\hat{\be}^{LS};\be)|
$$
if $1\leq\ell<(p-2)/2$.
\hfill$\Box$
\end{exm}

The following theorem specifies a general condition that $\hat{\be}_\phi$, given in \eqref{eqn:class}, reduces the bias of $\hat{\be}^{LS}$ in a finite sample setup.
\begin{thm}\label{thm:BR}
Assume that $p\geq 5$.
Let the $a_j$ and the $b_j$ be defined as in \eqref{eqn:BRi}.
Assume that $\phi(t)$ is bounded as $0\leq \phi(t)\leq 2 \sum_{j=1}^\ell (a_j/b_j)t^{-j}$ for any $t>0$ and a fixed natural number $\ell$.
If $\ell<(p-2)/2$, then we have $|{\rm Bias}(\hat{\be}_\phi;\be)|\leq |{\rm Bias}(\hat{\be}^{LS};\be)|$ for any $\be$.
\end{thm}

{\bf Proof.}\ \ 
Using the same arguments as in \eqref{eqn:EBR2}, we can express $|{\rm Bias}(\hat{\be}_\phi;\be)|$ as $|{\rm Bias}(\hat{\be}_\phi;\be)|=|-E_0+E_\phi|\cdot|\be|$, where
$$
E_0=E\bigg[\frac{p-2}{p+2K-2}\bigg],\quad E_\phi=E\bigg[\phi\bigg(\frac{\Vert\U\Vert^2}{S}\bigg)\frac{\U^t\bxi}{\Vert\U\Vert^2}\bigg].
$$
From Lemma \ref{lem:bias_LS}, it suffices to show that $|-E_0+E_\phi|\leq E_0$ or, equivalently, that
\begin{equation}\label{eqn:bias-phi1}
-2E_0\leq -2E_0+E_\phi\leq 0.
\end{equation}

Since $\phi(t)\geq 0$ for any $t$, it follows from (i) of Lemma \ref{lem:expectations} that $E_\phi\geq 0$.
Thus the first inequality of \eqref{eqn:bias-phi1} is valid.

Combining (i) of Lemma \ref{lem:expectations} and the given boundedness assumption on $\phi$ yields that
\begin{align*}
E_\phi
&=\sum_{k=0}^\infi \frac{2\la}{p+2k}P_\la(k)\int_0^\infi\!\!\!\int_0^\infi\phi\Big(\frac{w}{s}\Big)g_{p+2k}(w)\dd w \,g_m(s)\dd s \\
&\leq\sum_{k=0}^\infi \frac{2\la}{p+2k}P_\la(k)\int_0^\infi\!\!\!\int_0^\infi 2\bigg\{\sum_{j=1}^\ell\frac{a_j}{b_j}\Big(\frac{s}{w}\Big)^j\bigg\}g_{p+2k}(w)\dd w \,g_m(s)\dd s \\
&=2 E\bigg[\sum_{j=1}^\ell\frac{a_j}{b_j}\bigg(\frac{S}{\Vert\U\Vert^2}\bigg)^j\frac{\U^t\bxi}{\Vert\U\Vert^2}\bigg].
\end{align*}
Hence, by the same arguments as in the proof of Lemma \ref{lem:bias-BRi}, it is seen that
\begin{align*}
-2E_0+E_\phi&\leq -2E_0+2E\bigg[\sum_{j=1}^\ell\frac{a_j}{b_j}\bigg(\frac{S}{\Vert\U\Vert^2}\bigg)^j\frac{\U^t\bxi}{\Vert\U\Vert^2}\bigg]\\
&=-2E\bigg[\prod_{j=1}^{\ell+1}\frac{p-2j}{p+2K-2j}\bigg]\leq 0,
\end{align*}
which implies that the second inequality of \eqref{eqn:bias-phi1} is valid.
\qed

\begin{exm}
Let $\phi_\ell^{ST}(t)=\sum_{j=1}^\ell(p/m)^j t^{-j}$.
Estimator \eqref{eqn:ST} can be expressed as $\hat{\be}_\ell^{ST}=\{1+\phi_\ell^{ST}(\Vert\U\Vert^2/S)\}\hat\be^{LS}$.
It is observed that
\begin{align*}
2\frac{a_j}{b_j}-\Big(\frac{p}{m}\Big)^j
&=\frac{p}{m}\Big\{2\frac{p-2j}{p}\frac{m}{m+2j-2}\frac{a_{j-1}}{b_{j-1}}-\Big(\frac{p}{m}\Big)^{j-1}\Big\} \\
&\leq \frac{p}{m}\Big\{2\frac{a_{j-1}}{b_{j-1}}-\Big(\frac{p}{m}\Big)^{j-1}\Big\},
\end{align*}
which implies that, if $(p/m)^\ell\leq 2a_\ell/b_\ell$ for a given natural number $\ell <(p-2)/2$, it follows that $(p/m)^j\leq 2a_j/b_j$ for $j=1,\ldots,\ell-1$.
Hence, using Theorem \ref{thm:BR}, we can obtain $|{\rm Bias}(\hat{\be}_\ell^{ST};\be)|\leq |{\rm Bias}(\hat{\be}^{LS};\be)|$ if $(p/m)^\ell\leq 2a_\ell/b_\ell$ for a given natural number $\ell <(p-2)/2$.
\hfill$\Box$
\end{exm}

\begin{exm}
Denote
$$
\hat\be_{\ell\cdot 2}^{BR}=\bigg\{1+2\sum_{j=1}^\ell\frac{a_j}{b_j}\bigg(\frac{S}{\Vert\U\Vert^2}\bigg)^j\bigg\}\hat{\be}^{LS}.
$$
It holds that $|{\rm Bias}(\hat{\be}_{\ell\cdot 2}^{BR};\be)|\leq |{\rm Bias}(\hat{\be}^{LS};\be)|$.
However, the bias of $\hat\be_{\ell\cdot 2}^{BR}$ does not always have the same sign as that of $\hat\be^{LS}$.
\hfill$\Box$
\end{exm}

\begin{exm}\label{exm:imp_bias}
The first moment of $\hat{\be}^{MM}$ is not finite.
Such an estimator not having finite moments can be modified by Theorem \ref{thm:BR}.

Assume that an estimator of $\be$ has the form $\hat{\be}_{\bar{\phi}}=\{1+\bar{\phi}(\Vert\U\Vert^2/S)\}\hat{\be}^{LS}$.
Let $\hat{\be}_{\phi_\ell^*}=\{1+\phi_\ell^*(\Vert\U\Vert^2/S)\}\hat{\be}^{LS}$, where $\ell$ is a natural number and
$$
\phi_\ell^*(t)=\max\bigg[0,\min\bigg\{\bar{\phi}(t), \sum_{j=1}^\ell \frac{a_j}{b_j}t^{-j}\bigg\}\bigg].
$$
If $\ell<(p-2)/2$, then $\hat{\be}_{\phi_\ell^*}$ has a finite smaller bias than $\hat{\be}^{LS}$ for any $\be$.
\hfill$\Box$
\end{exm}

\begin{exm}
The second moment of $\hat{\be}_\ell^{BR}$ is always larger than that of $\hat{\be}^{LS}$.
Thus there is a considerable risk that $\hat{\be}_\ell^{BR}$ has larger variance and MSE than $\hat{\be}^{LS}$.
To reduce the risk, we consider, for example, the following truncation rule
$$
\phi_\ell^{**}(t)=\begin{cases}
\sum_{j=1}^\ell(a_j/b_j)t^{-j} & \textup{if $t>1$}, \\
(a_1/b_1)t^{-1} & \textup{otherwise}.
\end{cases}
$$
Then, the resulting estimator $\hat{\be}_{\phi_\ell^{**}}=\{1+\phi_\ell^{**}(\Vert\U\Vert^2/S)\}\hat{\be}^{LS}$ always has a smaller second moment than $\hat{\be}_\ell^{BR}$.
\hfill$\Box$
\end{exm}

\section{MSE reduction}\label{sec:MSE}

In estimation of a normal mean vector $\bxi$ with a quadratic loss, where $\U\sim\Nc_{p}(\bxi,\si^2I_p)$ and $S\sim\si^2\chi^2_m$, it is well known that the ML estimator, $\widehat{\bxi}{}^{ML}=\U$, is uniformly dominated by the James and Stein (1961) shrinkage estimator $\widehat{\bxi}{}^{JS}=(1-G^{JS})\U$ with $G^{JS}=(p-2)S/\{(m+2)\Vert\U\Vert^2\}$.
Moreover, from the integral expression of risk difference (IERD) method by Kubokawa (1994), we can show that $\widehat{\bxi}{}^{JS}$ is improved by a truncated shrinkage estimator $\widehat{\bxi}{}^{K}=(1-G^K)\U$ with $G^K=\min\{(p-2)/p,  G^{JS}\}$.

Whittemore (1989) and Guo and Ghosh (2012) employed the above shrinkage estimators to find out better slope estimators for a linear measurement error model with a structural relationship.
Their ideas can be applied to our slope estimation in the functional model \eqref{eqn:model2}.
For the ordinary LS estimator $\hat{\be}^{LS}=\U^t\Z/\Vert\U\Vert^2$, substituting $\U$ with $\widehat{\bxi}{}^{JS}$ yields Whittemore (1989) type estimator
$$
\hat{\be}^W
=\frac{(\widehat{\bxi}{}^{JS})^t\Z}{\Vert\widehat{\bxi}{}^{JS}\Vert^2}
=\frac{\U^t\Z}{(1-G^{JS})\Vert\U\Vert^2}.
$$
Similarly, by replacing $\U$ with $\widehat{\bxi}{}^{K}$, we obtain Guo and Ghosh (2012) type estimator
\begin{equation}\label{eqn:GG}
\hat{\be}^{GG}
=\frac{(\widehat{\bxi}{}^{K})^t\Z}{\Vert\widehat{\bxi}{}^{K}\Vert^2}
=\frac{\U^t\Z}{(1-G^{K})\Vert\U\Vert^2}.
\end{equation}
The Whittemore estimator $\hat{\be}^W$ is asymptotically analogous to the method of moments estimator $\hat{\be}^{MM}$ given in Section 3, and the bias and the MSE of $\hat{\be}^W$ do not exist.
Meanwhile, the Guo and Ghosh estimator $\hat{\be}^{GG}$ has a finite MSE.

In this section, a unified method is provided for the MSE reduction not only for $\hat{\be}^{LS}$ and $\hat{\be}^{GG}$, but also for the bias-reduced estimators $\hat{\be}_\phi$ given in Section \ref{sec:bias}.

\subsection{Preliminaries}

Suppose that an estimator of the slope $\be$ in reparametrized model \eqref{eqn:model2} depends only on $\Z$, $\U$ and $S$ but not on $Z_0$ and $U_0$.
Recall that
\begin{equation}\label{eqn:model3}
\Z\sim\Nc_p(\be\bxi,\tau^2I_p),\qquad
\U\sim\Nc_p(\bxi,\si^2I_p),\qquad
S\sim \si^2\chi_m^2.
\end{equation}
If $\tau^2=\si^2$ in partial model \eqref{eqn:model3}, the problem of estimating $\be$ is just the same as a linear calibration problem. 
More precisely, the MSE reduction problem for $\hat{\be}^{LS}$ corresponds to that for what is called a classical estimator in the multivariate linear calibration problem with a single independent variable.
For details of the linear calibration problem, see Kubokawa and Robert (1994), who derived an alternative to the classical estimator under the MSE criterion.
See also Osborne (1991), Brown (1993) and Sundberg (1999) for a general overview of the calibration problem.

Let $V=\Vert\U\Vert^2/(S+\Vert\U\Vert^2)$ and let $\psi(v)$ be a function on the interval $(0,1)$.
In this section, we consider an alternative estimator of the form
$$
\hat{\be}_\psi=\psi(V)\hat{\be}^{LS}=\psi(V)\frac{\U^t\Z}{\Vert\U\Vert^2}.
$$
It is clear that
\begin{align*}
{\rm MSE}(\hat{\be}_\psi;\be)&=E\bigg[\psi^2(V)\frac{\U^t\Z\Z^t\U}{\Vert\U\Vert^4}-2\be\psi(V)\frac{\U^t\Z}{\Vert\U\Vert^2}-\be^2 \bigg].
\end{align*}
Taking expectation with respect to $\Z\sim\Nc_p(\be\bxi,\tau^2I_p)$ gives that
\begin{equation}\label{eqn:MSEpsi}
{\rm MSE}(\hat{\be}_\psi;\be)=\tau^2E\bigg[\frac{\psi^2(V)}{\Vert\U\Vert^2}\bigg]+\be^2E\bigg[\bigg\{\psi(V)\frac{\U^t\bxi}{\Vert\U\Vert^2}-1\bigg\}^2\bigg].
\end{equation}
Hence, if $\psi^2(V)\leq 1$ and
\begin{equation}\label{eqn:p2}
E\bigg[\bigg\{\psi(V)\frac{\U^t\bxi}{\Vert\U\Vert^2}-1\bigg\}^2\bigg]\leq E\bigg[\bigg\{\frac{\U^t\bxi}{\Vert\U\Vert^2}-1\bigg\}^2\bigg],
\end{equation}
then $\hat{\be}_\psi$ has a smaller MSE than $\hat{\be}^{LS}$.

As pointed out by Kubokawa and Robert (1994), condition \eqref{eqn:p2} is closely related to a statistical control problem.
The control problem is formulated as the problem of estimating a normal mean vector $\bxi$, where the accuracy of an estimator $\hat{\bxi}$ is measured by loss $(\hat{\bxi}{}^t\bxi-1)^2$.
For more details of the statistical control problem, the reader is referred to Zellner (1971) and also to Zaman (1981), Berger et al. (1982) and Aoki (1989).

In Kubokawa and Robert (1994), the IERD method (Kubokawa (1994)) plays an important role in checking condition \eqref{eqn:p2}.
Here, we do not employ the IERD method and we directly evaluate the expectations in \eqref{eqn:p2} with the help of a Poisson variable.
\begin{lem}\label{lem:expect-control}
For nonnegative integers $k$, denote by $P_\la(k)$ the Poisson probabilities with mean $\la=\Vert\bxi\Vert^2/(2\si^2)$.
Assume that $\psi(V)\U^t\bxi/\Vert\U\Vert^2$ has a finite second moment.
Then we have
$$
E\bigg[\bigg\{\psi(V)\frac{\U^t\bxi}{\Vert\U\Vert^2}-1\bigg\}^2\bigg]-1
=\sum_{k=0}^\infi \frac{2\la}{p+2k}P_\la(k)H_\psi(k),
$$
where
\begin{align*}
H_\psi(k)&=\int_0^1\bigg\{\frac{1+2k}{p+2k+m-2}\frac{\psi^2(v)}{v}-2\psi(v)\bigg\}f_k(v)\dd v , \\
f_k(v)&=\frac{\Ga[(p+2k+m)/2]}{\Ga[(p+2k)/2]\Ga[m/2]}v^{(p+2k)/2-1}(1-v)^{m/2-1}.
\end{align*}
\end{lem}

{\bf Proof.}\ \ 
Note that $V$ can be interpreted as a function of $\Vert\U\Vert^2/S$.
For that reason, Lemma \ref{lem:expectations} can be used to obtain
$$
E\bigg[\bigg\{\psi(V)\frac{\U^t\bxi}{\Vert\U\Vert^2}-1\bigg\}^2\bigg]-1 
=\sum_{k=0}^\infi \frac{2\la}{p+2k}P_\la(k)H_\psi^*(k),
$$
where
$$
H_\psi^*(k)=\int_0^\infi\!\!\!\int_0^\infi\bigg\{\frac{1+2k}{w}\psi^2\Big(\frac{w}{w+s}\Big)-2\psi\Big(\frac{w}{w+s}\Big)\bigg\}g_m(s)g_{p+2k}(w)\dd s\dd w.
$$
For $H_\psi^*(k)$, we make the change of variables $t=s+w$ and $v=w/(w+s)$ with the Jacobian $J[(s,w)\to(t,v)]=t$ and hence
$$
H_\psi^*(k)=\int_0^1\!\!\int_0^\infi\bigg\{\frac{1+2k}{tv}\psi^2(v)-2\psi(v)\bigg\}g_{p+2k+m}(t)f_k(v)\dd t\dd v.
$$
Integrating out with respect to $t$ yields that $H_\psi^*(k)=H_\psi(k)$, which completes the proof.
\qed

\medskip
Next, we specify conditions for finiteness of the MSEs of $\hat\be^{LS}$ and $\hat\be_\ell^{BR}$, where $\hat\be_\ell^{BR}$ is given in \eqref{eqn:BRi}.
\begin{lem}\label{lem:MSE_LS}
Let $K$ be a Poisson random variable with mean $\la=\Vert\bxi\Vert^2/(2\si^2)$.
If $p\geq 3$, the MSE of $\hat\be^{LS}$ is finite and it can be expressed as
\begin{align}\label{eqn:MSE_LS}
{\rm MSE}(\hat{\be}^{LS};\be) 
&=\frac{\tau^2}{\si^2}E\bigg[\frac{1}{p+2K-2}\bigg]\non\\
&\qquad +\be^2E\bigg[\frac{2\la(1+2K)}{(p+2K)(p+2K-2)}-\frac{4\la}{p+2K}+1\bigg].
\end{align}
\end{lem}

{\bf Proof.}\ \ 
From \eqref{eqn:MSEpsi}, the MSE of $\hat\be^{LS}$ can be written as
$$
{\rm MSE}(\hat{\be}^{LS};\be)=\tau^2E\bigg[\frac{1}{\Vert\U\Vert^2}\bigg]+\be^2E\bigg[\frac{(\U^t\bxi)^2}{\Vert\U\Vert^4}-2\frac{\U^t\bxi}{\Vert\U\Vert^2}+1\bigg].
$$
Using identities \eqref{eqn:LS-1st} and \eqref{eqn:LS-2nd} and Lemma \ref{lem:chi}, we obtain the lemma.
\qed

\medskip
If $\la$ has the same order as $n$, the first term of the r.h.s. in \eqref{eqn:MSE_LS} converges to zero as $n\to\infi$.
Hence, the MSE of $\hat\be^{LS}$ is not much influenced by $\tau^2$ when $n$ is sufficiently large or when $\tau^2$ is sufficiently smaller than $\si^2$.

\begin{lem}\label{lem:MSE_BR}
Assume that $p\geq 7$.
If $1\leq \ell<(p-2)/4$, the MSE of $\hat\be_\ell^{BR}$ is finite.
\end{lem}

{\bf Proof.}\ \ 
From \eqref{eqn:MSEpsi}, it is sufficient to derive a condition that
$$
E\bigg[\bigg(\frac{S}{\Vert\U\Vert^2}\bigg)^{2\ell}\frac{1}{\Vert\U\Vert^2}\bigg]<\infi
\quad\textup{and}\quad
E\bigg[\bigg(\frac{S}{\Vert\U\Vert^2}\bigg)^{2\ell}\frac{(\U^t\bxi)^2}{\Vert\U\Vert^4}\bigg]<\infi.
$$

Lemma \ref{lem:chi} leads to, for $p-4\ell-2>0$,
\begin{align*}
E\bigg[\bigg(\frac{S}{\Vert\U\Vert^2}\bigg)^{2\ell}\frac{1}{\Vert\U\Vert^2}\bigg]
&=E\bigg[\frac{\si^{4\ell}}{\Vert\U\Vert^{4\ell+2}}\bigg] \prod_{j=1}^{2\ell}(m+2j-2)\\
&=\frac{1}{\si^2}E\bigg[\prod_{j=1}^{2\ell+1}\frac{1}{p+2K-2j}\bigg] \prod_{j=1}^{2\ell}(m+2j-2).
\end{align*}
Similarly, using (ii) of Lemma \ref{lem:expectations} and Lemma \ref{lem:chi} yields that for $p-4\ell-2>0$
$$
E\bigg[\bigg(\frac{S}{\Vert\U\Vert^2}\bigg)^{2\ell}\frac{(\U^t\bxi)^2}{\Vert\U\Vert^4}\bigg]
=E\bigg[\frac{2\la(1+2K)}{p+2K}\prod_{j=1}^{2\ell+1}\frac{1}{p+2K-2j}\bigg]\prod_{j=1}^{2\ell}(m+2j-2).
$$
Hence the finiteness of the MSE of $\hat\be_\ell^{BR}$ needs $p-4\ell-2>0$, namely $\ell<(p-2)/4$. 
\qed

\medskip
We can express the MSE of $\hat\be_\ell^{BR}$ alternatively by using the Poisson random variable as in Lemma \ref{lem:MSE_LS}, but it is omitted.

\subsection{Main analytical result and some examples}

Consider a slope estimator of the form $\hat{\be}_{\bar{\psi}}=\bar{\psi}(V)\hat{\be}^{LS}$, where $\bar{\psi}(v)$ is a function of $v$ on the interval $(0,1)$.
Assume that the second moment of $\hat{\be}_{\bar{\psi}}$ is finite.
Suppose that we want to find out an estimator $\hat{\be}_{\psi}=\psi(V)\hat{\be}^{LS}$ having a smaller MSE than $\hat{\be}_{\bar{\psi}}$, where $\psi(v)$ is a function on $(0,1)$.
To this end, $\psi$ requires some conditions in the following theorem.
\begin{thm}\label{thm:imp_be}
If $\psi^2(v)\leq \bar{\psi}^2(v)$ and 
$$
\De(v|\psi,\bar{\psi})=\psi^2(v)-\bar{\psi}^2(v)-2(p+m-2)v\{\psi(v)-\bar{\psi}(v)\}\leq 0
$$
for any $v\in (0,1)$, then ${\rm MSE}(\hat{\be}_{\psi};\be)\leq {\rm MSE}(\hat{\be}_{\bar{\psi}};\be)$.
\end{thm}

\noindent
{\bf Proof.}\ \ 
Since $\psi^2(v)\leq \bar{\psi}^2(v)$ for any $v$, $\hat\be_\psi$ inherits the finiteness of the second moment from $\hat{\be}_{\bar{\psi}}$.
By virtue of Lemma \ref{lem:expect-control}, the difference between the MSEs of $\hat{\be}_{\psi}$ and $\hat{\be}_{\bar{\psi}}$ is expressed as
\begin{align*}
{\rm MSE}(\hat{\be}_{\psi};\be)-{\rm MSE}(\hat{\be}_{\bar{\psi}};\be)
&=\tau^2E\bigg[\frac{\psi^2(V)-\bar{\psi}^2(V)}{\Vert\U\Vert^2}\bigg] \\
&\qquad +\be^2\sum_{k=0}^\infi \frac{2\la}{p+2k}P_\la(k)\int_0^1\De_k(v|\psi,\bar{\psi})\frac{f_k(v)}{v}\dd v,
\end{align*}
where
$$
\De_k(v|\psi,\bar{\psi})=\frac{1+2k}{p+2k+m-2}\{\psi^2(v)-\bar{\psi}^2(v)\}-2v\{\psi(v)-\bar{\psi}(v)\}.
$$
It follows that for any $k\geq 0$
$$
\frac{1+2k}{p+2k+m-2}\geq \frac{1}{p+m-2},
$$
which implies that $\De_k(v|\psi,\bar{\psi})\leq \De_0(v|\psi,\bar{\psi})=\De(v|\psi,\bar{\psi})/(p+m-2)$.
Hence the proof is complete.
\qed

\medskip
Theorem \ref{thm:imp_be} is the key to constructing a better estimator under the MSE criterion.
In the following, we show some examples.
\begin{exm}\label{exm:imp_be0}
For a given $\bar{\psi}$, let $\psi_0(v)=\max[0,\min\{\bar{\psi}(v), 2(p+m-2)v-\bar{\psi}(v)\}]$.
Note that $0\leq\psi_0(v)\leq|\bar{\psi}(v)|$ for any $v\in (0,1)$.
It also turns out that
$$
\psi_0(v)=\left\{\begin{array}{ll}
2(p+m-2)v-\bar{\psi}(v) & \begin{array}{l} \textup{if $0\leq 2(p+m-2)v-\bar{\psi}(v)\leq \bar{\psi}(v)$}, \end{array} \\
\bar{\psi}(v) & \begin{array}{l}\textup{if $0\leq \bar{\psi}(v) < 2(p+m-2)v-\bar{\psi}(v)$}, \end{array} \\
0 & \begin{array}{l}
       \textup{if $2(p+m-2)v-\bar{\psi}(v)< 0\leq \bar{\psi}(v)$} \\
       \textup{\quad or if $\bar{\psi}(v)<0\leq 2(p+m-2)v-\bar{\psi}(v)$},
     \end{array}
\end{array}
\right.
$$
which implies that $\De(v|\psi_0,\bar{\psi})=0$ if $\min\{\bar{\psi}(v), 2(p+m-2)v-\bar{\psi}(v)\}\geq 0$ and $\De(v|\psi_0,\bar{\psi})=\bar{\psi}(v)\{2(p+m-2)v-\bar{\psi}(v)\}\leq 0$ otherwise.
Hence, if $\Pr(\psi_0(V)=\bar{\psi}(V))<1$, then $\hat{\be}_{\psi_0}=\psi_0(V)\hat\be^{LS}$ is better than $\hat{\be}_{\bar{\psi}}$ under the MSE criterion.

Particularly, when $\bar\psi(v)\equiv 1$,
\begin{equation}\label{eqn:TLS}
\hat{\be}^{TLS}=\max\bigg[0,\min\bigg\{1,\, \frac{2(p+m-2)\Vert\U\Vert^2}{S+\Vert\U\Vert^2}-1\bigg\}\bigg]\hat{\be}^{LS}
\end{equation}
has a smaller MSE than $\hat{\be}^{LS}$ for $p\geq 3$.
\hfill$\Box$
\end{exm}

\begin{exm}\label{exm:imp_be1}
Assume additionally that $\bar{\psi}(v)\geq 1$ for $0<v<1$.
Let
$$
\psi_1(v)=\max[1,\min\{\bar{\psi}(v), 2(p+m-2)v-\bar{\psi}(v)\}].
$$
Using the same arguments as in Example \ref{exm:imp_be0}, we can prove that if $\Pr(\psi_1(V)=\bar{\psi}(V))<1$ then $\hat{\be}_{\psi_1}$ has a smaller MSE than $\hat{\be}_{\bar{\psi}}$.
Since $\psi_1(v)\geq 1$ for any $v\in(0,1)$, it holds true that $|E[\hat{\be}_{\psi_1}]|\geq|E[\hat\be^{LS}]|$, which implies that $\hat{\be}_{\psi_1}$ not only improves on the MSE of $\hat{\be}_{\bar{\psi}}$, but also may correct the bias of $\hat\be^{LS}$.
\hfill$\Box$
\end{exm}

\begin{exm}\label{exm:imp_be2}
Assume that $p\geq 7$.
Let 
$$
\bar\psi_\ell^{BR}(v)=1+\sum_{j=1}^\ell \frac{a_j}{b_j}\bigg(\frac{1-v}{v}\bigg)^j.
$$
The MSE of the bias-reduced estimator $\hat\be_\ell^{BR}=\bar\psi_\ell^{BR}(V)\hat\be^{LS}$ is improved by
\begin{equation}\label{eqn:TBRi}
\begin{split}
\hat{\be}_\ell^{TBR}&=\psi_\ell^{TBR}(V)\hat{\be}^{LS},\\
\psi_\ell^{TBR}(v)&=\max[1,\min\{\bar{\psi}_\ell^{BR}(v),2(p+m-2)v-\bar{\psi}_\ell^{BR}(v)\}]
\end{split}
\end{equation}
for $1\leq\ell<(p-2)/4$.

Guo and Ghosh's (2012) estimator can be written as $\hat{\be}^{GG}=\bar{\psi}^{GG}(V)\hat\be^{LS}$ with
$$
\bar{\psi}^{GG}(v)=\bigg[1-\min\bigg\{\frac{p-2}{p},\frac{p-2}{m+2}\frac{1-v}{v}\bigg\}\bigg]^{-1}.
$$
Define
\begin{equation}\label{eqn:TGG}
\begin{split}
\hat{\be}^{TGG}&=\psi^{TGG}(V)\hat{\be}^{LS}, \\
\psi^{TGG}(v)&=\max[1,\min\{\bar{\psi}^{GG}(v),2(p+m-2)v-\bar{\psi}^{GG}(v)\}].
\end{split}
\end{equation}
Since $\Pr(\psi^{TGG}(V)=\bar{\psi}^{GG}(V))< 1$, $\hat{\be}^{TGG}$ dominates $\hat{\be}^{GG}$ under the MSE criterion.
\hfill$\Box$
\end{exm}

\begin{exm}\label{exm:imp_be3}
An improved estimator on $\hat{\be}^{LS}$ can be obtained by means of Equation (2.4) of Kubokawa and Robert (1994).

Assume that $\bar\psi(v)\geq0$ for $0<v<1$.
Let $\psi^{KR}(v)=\min\{\bar\psi(v), (p+m-2)v\}$.
Then it is easy to show from Theorem \ref{thm:imp_be} that $\hat{\be}^{KR}=\psi^{KR}(V)\hat{\be}^{LS}$ has a smaller MSE than $\hat{\be}_{\bar\psi}$ when $\Pr(\psi^{KR}(V)=\bar\psi(V))<1$.
From the above-mentioned, it is obvious that
$$
\hat{\be}^{TLS2}
=\min\{1,(p+m-2)V\}\hat{\be}^{LS}
=\min\bigg\{\frac{1}{\Vert\U\Vert^2},\,\frac{p+m-2}{\Vert\U\Vert^2+S}\bigg\}\U^t\Z
$$
has a smaller MSE than $\hat{\be}^{LS}$ for $p\geq 3$.
The estimator $\hat{\be}^{TLS2}$ is quite similar to an estimator given in Corollary 2.2 of Kubokawa and Robert (1994).
\hfill$\Box$
\end{exm}

\section{Numerical studies}\label{sec:num}
\subsection{Numerical examples with corn yield data}

In this subsection, numerical examples with real data sets illustrate how regression lines are drawn with the LS and its bias-reduced estimates and also with the ML and the inverse regression estimates.

For simplicity, we suppose $\tau^2=\si_x^2\ (=r\si^2)$ in model \eqref{eqn:model0}.
Then, the ML estimator of $\be$ has the form
\begin{equation}\label{eqn:ML}
\hat{\be}^{ML}=\frac{\Vert\Z\Vert^2-r\Vert\U\Vert^2+\sqrt{(\Vert\Z\Vert^2-r\Vert\U\Vert^2)^2+4r(\U^t\Z)^2}}{2\U^t\Z},
\end{equation}
which can be constructed by minimizing
$$
\frac{1}{\tau^2}\Vert\Z-\be\bxi\Vert^2+\frac{1}{\si^2}\Vert\U-\bxi\Vert^2
=\frac{1}{\tau^2}\{\Vert\Z-\be\bxi\Vert^2+r\Vert\U-\bxi\Vert^2\}
$$
subject to $-\infi<\be<\infi$ and $\bxi\in\Rb^p$.
Under a suitable convergence condition, $\hat{\be}^{ML}$ is a consistent estimator of $\be$.

As stated in the beginning of Subsection \ref{subsec:reparametrize}, $\hat\be^{LS}$ is derived from the regression of the $Y_i$ on the $\overline{X}_i$.
Let us now consider the inverse regression, namely the $\overline{X}_i$ are regressed on the $Y_i$.
Through the use of statistics in \eqref{eqn:model2}, the least squares estimator for a slope of the inverse regression equals to $\U^t\Z/\Vert\Z\Vert^2$.
Since the slope of the inverse regression is equivalent to $\be^{-1}$ (the reciprocal of the slope in the usual regression), the resulting estimator of $\be$ can be expressed as
\begin{equation}\label{eqn:IR}
\hat{\be}^{IR}=\frac{\Vert\Z\Vert^2}{\U^t\Z}.
\end{equation}

Note that $\hat{\be}^{ML}$ and $\hat{\be}^{IR}$ have no finite moments, and hence their biases and MSEs do not exist.
If $\U^t\Z>0$, then it can easily be shown that $(\hat{\be}^{ML})^{-1}<(\hat{\be}^{LS})^{-1}$ and $\hat{\be}^{ML}<\hat{\be}^{IR}$, namely  $0<\hat{\be}^{LS}<\hat{\be}^{ML}<\hat{\be}^{IR}$.
In a similar fashion, we obtain $\hat{\be}^{IR}<\hat{\be}^{ML}<\hat{\be}^{LS}<0$ if $\U^t\Z<0$.
See Anderson (1976).

\begin{figure}[t]
\begin{center}
\begin{minipage}[t]{.47\textwidth}
\begin{center}
\includegraphics[width=1.0\textwidth]{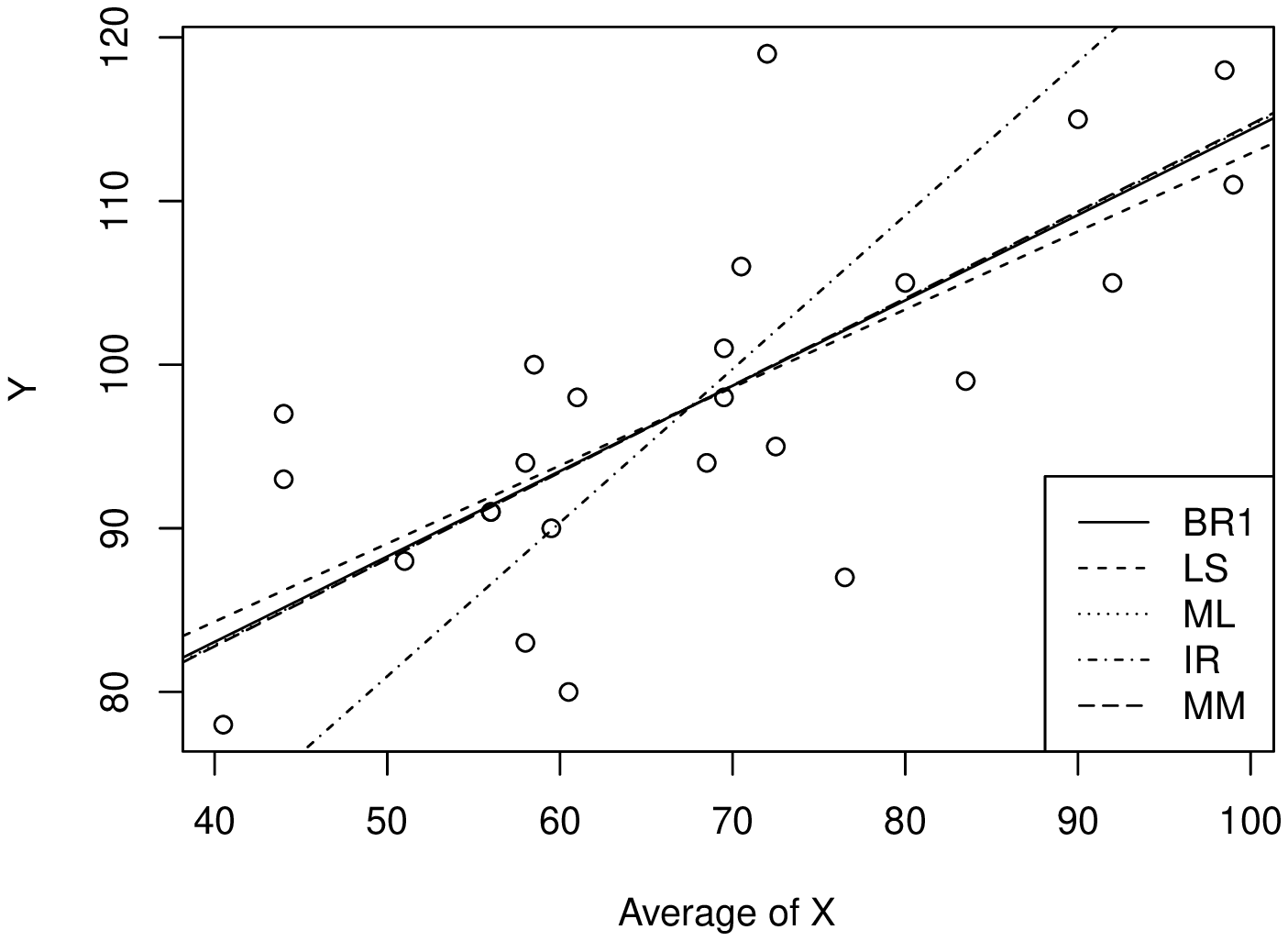}
\vspace{-30pt}
\caption{Some regression lines for Fuller's (1987) corn-yield data, where $n=25$ and $r=2$.}
\label{fig:1}
\end{center}
\end{minipage}
\hfill
\begin{minipage}[t]{.47\textwidth}
\begin{center}
\includegraphics[width=1.0\textwidth]{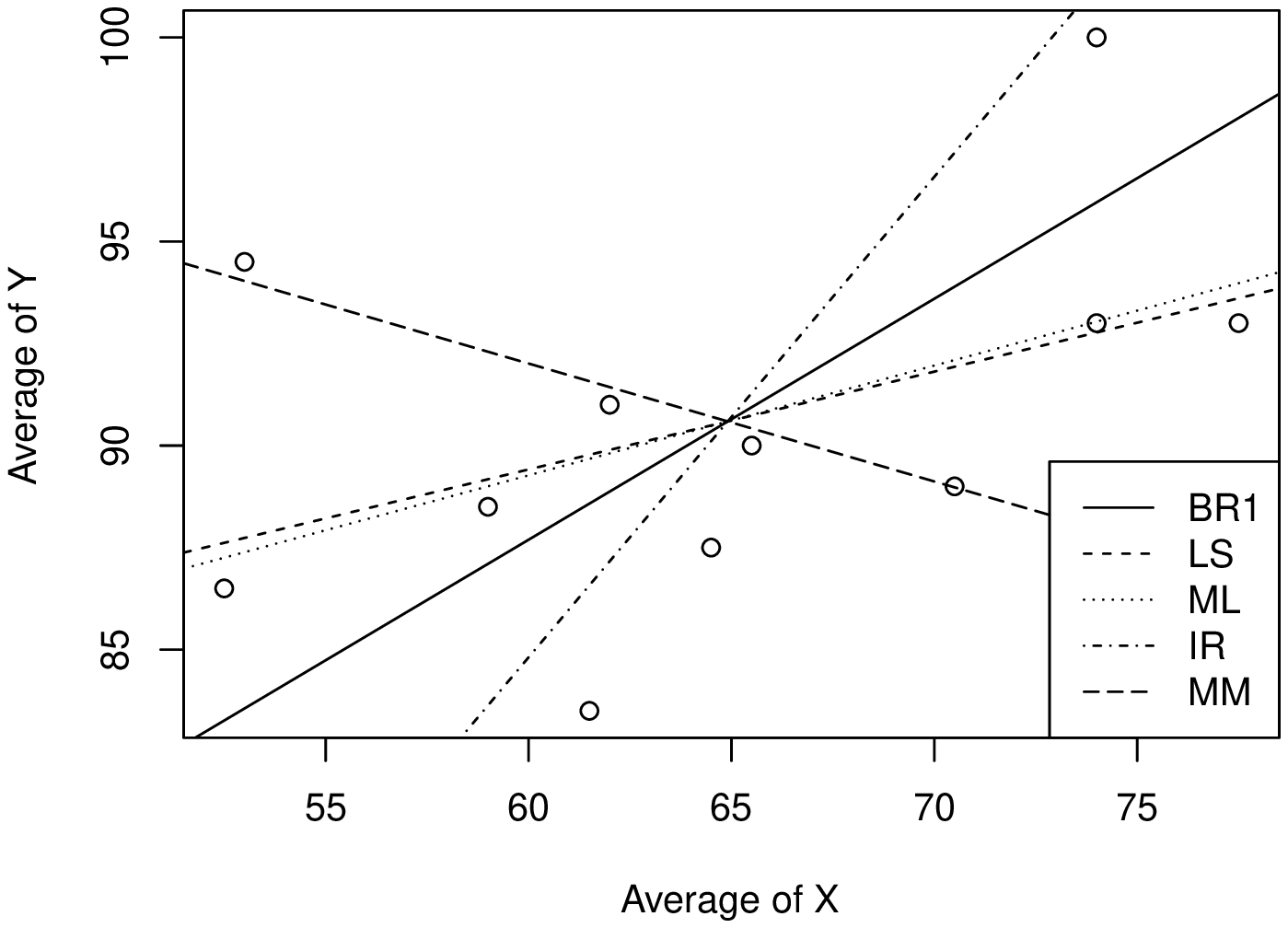}
\vspace{-30pt}
\caption{Some regression lines for DeGracie and Fuller's (1972) corn-yield data, where $n=11$ and $r=2$.}
\label{fig:2}
\end{center}
\end{minipage}
\end{center}
\end{figure}

We now present two numerical examples for corn-yield data sets given in Fuller (1987, Table 3.1.1) and in DeGracie and Fuller (1972, p.934). 
The data sets consist of the yields of corn with two soil nitrogen contents.
The yield and the soil nitrogen content are assumed to be, respectively, dependent $(Y)$ and independent $(X)$ variables, where the data set of DeGracie and Fuller (1972) has duplicate observations of the yield and so the average of the two yields was regarded as one dependent variable.
Figures \ref{fig:1} and \ref{fig:2} are scatter plots of the two data sets.
In the figures, we added several regression lines by using the ordinary LS estimate $(\hat\al^{LS},\hat\be^{LS})$, the bias-reduced (BR1) estimate $(\hat\al_1^{BR},\hat\be_1^{BR})$, the ML estimate $(\hat\al^{ML},\hat\be^{ML})$, the inverse regression (IR) estimate $(\hat\al^{IR},\hat\be^{IR})$ and the method of moments (MM) estimate $(\hat\al^{MM},\hat\be^{MM})$, where $\hat\be_1^{BR}$, $\hat\be^{LS}$, $\hat\be^{ML}$, $\hat\be^{IR}$ and $\hat\be^{MM}$ are defined in \eqref{eqn:LS}, \eqref{eqn:BRi}, \eqref{eqn:ML}, \eqref{eqn:IR} and \eqref{eqn:MM}, respectively, and the corresponding estimates of $\al$ in any procedure are $Z_0-\hat\be U_0$.
Also, Tables \ref{tab:1} and \ref{tab:2} give the above estimates for the two data sets and show how $\hat\be_\ell^{BR}$ changes as $\ell$ increases, where, in the tables, ${\rm BR}_\ell$ denotes $\hat\be_\ell^{BR}$.

\begin{table}[t]
\begin{center}
\begin{minipage}[t]{.48\textwidth}
\begin{center}
\caption{Estimates of the slope and the intercept parameters for Fuller's (1987) corn-yield data, where $n=25$ and $r=2$.}
\label{tab:1}
\vspace{6pt}
$
\begin{array}{lrr}
\hline
{\rm Procedure} & \multicolumn{1}{c}{\be} & \multicolumn{1}{c}{\al} \\
\hline
{\rm LS}   & 0.47693 & 65.219 \\[2pt]
{\rm BR}_1 & 0.52183 & 62.185 \\
{\rm BR}_2 & 0.52587 & 61.912 \\
{\rm BR}_3 & 0.52623 & 61.888 \\
{\rm BR}_4 & 0.52626 & 61.886 \\
{\rm BR}_5 & 0.52627 & 61.885 \\
{\rm BR}_6 & 0.52627 & 61.885 \\[2pt]
{\rm ML}   & 0.52860 & 61.728 \\
{\rm IR}   & 0.93854 & 34.032 \\
{\rm MM}   & 0.53151 & 61.531 \\
\hline
\end{array}
$
\end{center}
\end{minipage}
\hfill
\begin{minipage}[t]{.48\textwidth}
\begin{center}
\caption{Estimates of the slope and the intercept parameters for DeGracie and Fuller's (1972) corn-yield data, where $n=11$ and $r=2$.}
\label{tab:2}
\vspace{6pt}
$
\begin{array}{lrr}
\hline
{\rm Procedure} & \multicolumn{1}{c}{\be} & \multicolumn{1}{c}{\al} \\
\hline
{\rm LS}   &  0.23972  &  75.031 \\[2pt]
{\rm BR}_1 &  0.59055  &  52.259 \\
{\rm BR}_2 &  1.03857  &  23.179 \\
{\rm BR}_3 &  1.55946  & -10.632 \\[2pt]
{\rm ML}   &  0.26902  &  73.129 \\
{\rm IR}   &  1.17756  &  14.157 \\
{\rm MM}   & -0.28904  & 109.353 \\
\hline
\end{array}
$
\end{center}
\end{minipage}
\end{center}
\end{table}

Table \ref{tab:1} and Figure \ref{fig:1} indicate that $\hat{\be}_1^{BR}$, $\hat{\be}^{ML}$ and $\hat{\be}^{MM}$ take similar values, while Table \ref{tab:2} and Figure \ref{fig:2} show that they are very different.
Even though it theoretically follows that $0<\hat{\be}^{LS}<\hat{\be}^{ML}$ for $\U^t\Z>0$, ${\rm ML}$ is just slightly larger than ${\rm LS}$ as for the two data sets.

From the data set of DeGracie and Fuller (1972), the value of $S/\Vert\U\Vert^2$ is approximately $1421.5/706.41\approx 2$ and, as in Table \ref{tab:2}, $\hat\be_1^{BR}$ is calculated as
\begin{align*}
\hat\be_1^{BR}&=\Big(1+\frac{p-2}{m}\frac{S}{\Vert\U\Vert^2}\Big)\hat\be^{LS} \\
&=
\Big(1+\frac{10-2}{11}\times \frac{1421.5}{706.41}\Big)\times 0.23972=0.59055.
\end{align*}
When $S/\Vert\U\Vert^2$ takes a large value, the value of ${\rm BR}_\ell$ increases or decreases progressively as $\ell$ increases, and only the method of moments estimate for the slope has different sign from other estimates.
Furthermore the value of slope estimate impacts an intercept estimate as long as we use $Z_0-\hat\be U_0$ as an estimate of the intercept.

\subsection{Monte Carlo studies for bias and MSE comparison}

Next, some results of Monte Carlo simulations are provided in order to compare the biases and MSEs of slope estimators.

For three different sample sizes $n=10$, $30$ and $100$ with $r=2$, each of the simulated biases and MSEs is based on $500,000$ independent replications of $(\Z,\U,S)$.
It was assumed that $\be=-5$, $\tau^2=10$, and $\si^2=1$ or $10$.
For the latent variable $\bxi$, all the elements of $\bxi$ were set to be $1/\sqrt{10}$ or $\sqrt{5}$, namely $\Vert\bxi\Vert^2=p/10$ or $5p$, which implies that $\si_\xi^2\equiv\lim_{n\to\infi} \Vert\bxi\Vert^2/p=1/10$ or $5$.

Table \ref{tab:3} shows some values of $\la=\Vert\bxi\Vert^2/(2\si^2)$ which were assumed for our simulation.
For example, the smallest value of $\la$ is $0.045$ when $n=10$, $\si_\xi^2=1/10$ and $\si^2=10$, and the largest value of $\la$ is $247.5$ when $n=100$, $\si_\xi^2=5$ and $\si^2=1$.
\begin{table}[t]
\begin{center}
\caption{Mean $\la=\Vert\bxi\Vert^2/(2\si^2)$ of the Poisson distribution ($p=n-1$).}
\label{tab:3}
\vspace{6pt}
$
\begin{array}{cr@{\hspace{20pt}}D{.}{.}{-1}D{.}{.}{-1}D{.}{.}{-1}}
\hline
&& \multicolumn{3}{c}{n}\\
\cline{3-5}
\si_\xi^2 & \si^2 & \multicolumn{1}{c}{10} & \multicolumn{1}{c}{30} & \multicolumn{1}{c}{100} \\
\hline
1/10& 1&  0.45 &  1.45 &   4.95 \\
1/10&10&  0.045&  0.145&   0.495\\
5   & 1& 22.5  & 72.5  & 247.5  \\
5   &10&  2.25 & 7.25  &  24.75 \\
\hline
\end{array}
$
\end{center}
\end{table}

Slope estimators which were investigated in our simulation are $\hat{\be}^{LS}$, $\hat{\be}^{TLS}$, $\hat{\be}_\ell^{BR}$ $(\ell=1,\ 5)$, $\hat{\be}_\ell^{TBR}$ $(\ell=1,\ 5)$, $\hat{\be}^{GG}$ and $\hat{\be}^{TGG}$, which are given in \eqref{eqn:LS}, \eqref{eqn:TLS}, \eqref{eqn:BRi}, \eqref{eqn:TBRi}, \eqref{eqn:GG} and \eqref{eqn:TGG}, respectively.
The simulated biases and MSEs of the above estimators are summarized in Table \ref{tab:4}, where ${\rm LS}$, ${\rm TLS}$, ${\rm BR}_\ell$ $(\ell=1,\ 5)$, ${\rm TBR}_\ell$ $(\ell=1,\ 5)$, ${\rm GG}$ and ${\rm TGG}$ denote, respectively, $\hat{\be}^{LS}$, $\hat{\be}^{TLS}$, $\hat{\be}_\ell^{BR}$ $(\ell=1,\ 5)$, $\hat{\be}_\ell^{TBR}$ $(\ell=1,\ 5)$, $\hat{\be}^{GG}$ and $\hat{\be}^{TGG}$.
Since ${\rm BR}_5$ and ${\rm TBR}_5$ have no finite moments for $n=10$, we omitted them from our simulation.

\begin{table}[ht]
\begin{center}
\caption{Simulated bias and MSE in slope estimation for $\be=-5$.}
\label{tab:4}
{\small $
\begin{array}{cclllrrcrrcrr}
\hline
&&&&&\multicolumn{2}{c}{n=10} && \multicolumn{2}{c}{n=30} && \multicolumn{2}{c}{n=100} \\
\cline{6-7}
\cline{9-10}
\cline{12-13}
\si_\xi^2 & \si^2 &\multicolumn{2}{c}{{\rm Estimator}} 
 && \multicolumn{1}{c}{{\rm Bias}} & \multicolumn{1}{c}{{\rm MSE}}
 && \multicolumn{1}{c}{{\rm Bias}} & \multicolumn{1}{c}{{\rm MSE}}
 && \multicolumn{1}{c}{{\rm Bias}} & \multicolumn{1}{c}{{\rm MSE}} \\
\hline
1/10& 1&&{\rm LS}    && 4.54 & 22.17 &&  4.54 &  21.03 &&  4.54 &  20.76 \\
    &  &&{\rm TLS}   && 4.54 & 22.16 &&  4.54 &  21.03 &&  4.54 &  20.76 \\[2pt]
    &  &&{\rm BR}_1  && 4.12 & 27.65 &&  4.12 &  18.67 &&  4.13 &  17.48 \\
    &  &&{\rm TBR}_1 && 4.16 & 24.00 &&  4.12 &  18.67 &&  4.13 &  17.48 \\[2pt]
    &  &&{\rm BR}_5  &&      &       &&  2.77 & 108.32 &&  2.81 &  12.85 \\
    &  &&{\rm TBR}_5 &&      &       &&  3.08 &  24.68 &&  2.81 &  12.84 \\[2pt]
    &  &&{\rm GG}    && 3.66 & 33.88 &&  1.52 &  49.46 && -4.10 & 197.64 \\
    &  &&{\rm TGG}   && 3.72 & 31.06 &&  1.52 &  49.38 && -4.10 & 197.64 \\
\hline
1/10&10&&{\rm LS}    && 4.95 & 24.69 &&  4.95 &  24.55 &&  4.95 &  24.52 \\
    &  &&{\rm TLS}   && 4.95 & 24.69 &&  4.95 &  24.55 &&  4.95 &  24.52 \\[2pt]
    &  &&{\rm BR}_1  && 4.90 & 25.41 &&  4.90 &  24.22 &&  4.90 &  24.07 \\
    &  &&{\rm TBR}_1 && 4.91 & 24.87 &&  4.90 &  24.22 &&  4.90 &  24.07 \\[2pt]
    &  &&{\rm BR}_5  &&      &       &&  4.71 &  49.35 &&  4.71 &  22.86 \\
    &  &&{\rm TBR}_5 &&      &       &&  4.76 &  24.59 &&  4.71 &  22.85 \\[2pt]
    &  &&{\rm GG}    && 4.85 & 25.94 &&  4.57 &  26.39 &&  3.64 &  30.47 \\
    &  &&{\rm TGG}   && 4.86 & 25.57 &&  4.57 &  26.37 &&  3.64 &  30.47 \\
\hline
5   & 1&&{\rm LS}    && 0.70 &  1.02 &&  0.79 &   0.78 &&  0.82 &   0.72 \\
    &  &&{\rm TLS}   && 0.70 &  1.02 &&  0.79 &   0.78 &&  0.82 &   0.72 \\[2pt]
    &  &&{\rm BR}_1  && 0.08 &  1.09 &&  0.12 &   0.33 &&  0.13 &   0.11 \\
    &  &&{\rm TBR}_1 && 0.08 &  1.09 &&  0.12 &   0.33 &&  0.13 &   0.11 \\[2pt]
    &  &&{\rm BR}_5  &&      &       &&  0.00 &   0.39 &&  0.00 &   0.12 \\
    &  &&{\rm TBR}_5 &&      &       &&  0.00 &   0.39 &&  0.00 &   0.12 \\[2pt]
    &  &&{\rm GG}    && 0.08 &  1.29 &&  0.04 &   0.39 &&  0.01 &   0.11 \\
    &  &&{\rm TGG}   && 0.08 &  1.29 &&  0.04 &   0.39 &&  0.01 &   0.11 \\
\hline
5   &10&&{\rm LS}    && 3.25 & 11.22 &&  3.31 &  11.09 &&  3.33 &  11.10 \\
    &  &&{\rm TLS}   && 3.25 & 11.22 &&  3.31 &  11.09 &&  3.33 &  11.10 \\[2pt]
    &  &&{\rm BR}_1  && 2.06 &  8.75 &&  2.17 &   5.46 &&  2.21 &   5.06 \\
    &  &&{\rm TBR}_1 && 2.11 &  7.75 &&  2.17 &   5.46 &&  2.21 &   5.06 \\[2pt]
    &  &&{\rm BR}_5  &&      &       &&  0.37 &  20.98 &&  0.42 &   2.81 \\
    &  &&{\rm TBR}_5 &&      &       &&  0.44 &  11.20 &&  0.42 &   2.81 \\[2pt]
    &  &&{\rm GG}    && 0.80 & 13.57 && -2.05 &  62.49 && -2.15 & 132.24 \\
    &  &&{\rm TGG}   && 0.89 & 12.45 && -2.05 &  62.48 && -2.15 & 132.24 \\
\hline
\end{array}
$}
\end{center}
\end{table}

Lemma \ref{lem:bias-BRi} suggests that the bias of $\hat\be_\ell^{BR}$ is small for a large $\la$.
This has been confirmed by our simulations.
In particular, when $\la$ is large, ${\rm BR}_1$ and ${\rm BR}_5$ substantially improve not only the bias of LS but also its MSE.
When $\la$ is very small ($n=10$, $\si_\xi^2=1/10$ and $\si^2=10$), ${\rm BR}_1$ slightly improves on the bias of LS, while the MSE of ${\rm BR}_1$ is larger than that of ${\rm LS}$.
Also, as $n$ increases, the MSEs of ${\rm BR}_1$ and ${\rm BR}_5$ decrease and their absolute values of biases gradually increase, which implies that the variances of ${\rm BR}_1$ and ${\rm BR}_5$ decrease with increasing $n$.

${\rm TLS}$ causes only very slight decrease in MSE of ${\rm LS}$.
On the other hand, ${\rm TBR}_5$ makes successful reduction in MSE of ${\rm BR}_5$ and, particularly, the reduction is substantial when $n=30$.
This suggests that the truncation rule \eqref{eqn:TBRi} is notably effective in a higher-order bias-reduced estimator.

When $n=10$, ${\rm TGG}$ makes the MSE improvement on ${\rm GG}$ at the cost of bias.
Only ${\rm GG}$ and ${\rm TGG}$ have underestimated $\be$ in some cases.
Although ${\rm GG}$ has the MSE convergence to $\be$ under a structural model, the convergence rate is probably just a bit low.

\section{Remarks}\label{sec:remarks}

This paper considered a simple linear regression model with measurement error and discussed the bias and MSE reduction for slope estimation in a finite sample situation.
We conclude this paper with some remarks.
\newcounter{remark}
\begin{list}{}{
\topsep=4pt
\parsep=0pt
\parskip=0pt
\itemsep=4pt
\itemindent=0pt
\labelwidth=20pt
\labelsep=0pt
\leftmargin=20pt
\listparindent=6pt
\usecounter{remark}
}
\renewcommand{\makelabel}{%
(\roman{remark})
}
\item
For the simple linear regression model \eqref{eqn:model0}, we assume that $\si_x^2$ is known.
Then, it is assumed that $\si^2=\si_x^2/r=1$ without loss of generality, and model \eqref{eqn:model0} can be reduced to
\begin{equation}\label{eqn:model4}
\begin{aligned}
Z_0&\sim\Nc(\al+\be\th,\tau^2),\\
U_0&\sim\Nc(\th,1),
\end{aligned}
\qquad
\begin{aligned}
\Z&\sim\Nc_p(\be\bxi,\tau^2I_p),\\
\U&\sim\Nc_p(\bxi,I_p),
\end{aligned}
\end{equation}
where $Z_0$, $\Z$, $U_0$ and $\U$ are mutually independent, and $\al$, $\be$, $\th$, $\tau^2$ and $\bxi$ are unknown parameters.
For such a known-$\si_x^2$ case, we can use the same arguments as in Sections \ref{sec:bias} and \ref{sec:MSE} to improve on the bias or the MSE of an ordinary LS estimator even if $r=1$.
For further detail, see Appendix \ref{sec:A}.

\item
Consider here a simple structural model, where the latent variables $\th$ and $\bxi$ follow certain specified probability distributions.
Then reparametrized model \eqref{eqn:model2} is replaced with a conditional model:
\begin{equation*}
\begin{aligned}
Z_0|\th&\sim\Nc(\al+\be\th,\tau^2),\\
U_0|\th&\sim\Nc(\th,\si^2),
\end{aligned}
\qquad
\begin{aligned}
\Z|\bxi&\sim\Nc_p(\be\bxi,\tau^2I_p),\\
\U|\bxi&\sim\Nc_p(\bxi,\si^2I_p),
\end{aligned}
\qquad
\begin{aligned}
&\\
S&\sim \si^2\chi_m^2.
\end{aligned}
\end{equation*}
They are conditionally independent given $\th$ and $\bxi$.
Let $\hat{\be}_s^{LS}=\U^t\Z/\Vert\U\Vert^2$, which is the ordinary LS estimator of $\be$.
Denote by $E[\cdot|\th,\bxi]$ a conditional expectation with respect to $(Z_0,\Z,U_0,\U,S)$ given $\th$ and $\bxi$ and by $E^{\th,\bxi}[\cdot]$ an expectation with respect to $\th$ and $\bxi$.
The bias and MSE of $\hat{\be}_s^{LS}$ can be written, respectively, as
\begin{align*}
{\rm Bias}(\hat{\be}_s^{LS};\be) &=E^{\th,\bxi}[E[\hat{\be}_s^{LS}-\be|\th,\bxi]],\\
{\rm MSE}(\hat{\be}_s^{LS};\be) &=E^{\th,\bxi}[E[(\hat{\be}_s^{LS}-\be)^2|\th,\bxi]].
\end{align*}
Hence, it is possible to analytically improve the bias or the MSE of $\hat{\be}_s^{LS}$ by means of the reducing methods considered in this paper.

\item
In this paper, the MSE reduction of an estimator is based on shrinking the estimator toward zero, while the bias reduction is achieved by expanding the estimator.
A theoretically exact result on simultaneous reduction for both bias and MSE is still not known in a finite sample situation.

\item
Estimation of the intercept $\al$ in reparametrized model \eqref{eqn:model2} is an interesting problem.
Using the same arguments as in Sections \ref{sec:bias} and \ref{sec:MSE}, we can easily make the bias and MSE reduction of the LS estimator.

Define a class of estimators for $\al$ as $\hat{\al}_\phi=Z_0-\hat{\be}_\phi U_0$, where $\hat{\be}_\phi$ is given in \eqref{eqn:class}.
Note that $\hat{\be}_\phi$ is independent of $Z_0$ and $U_0$.
The bias of $\hat{\al}_\phi$ is written as
\begin{align*}
{\rm Bias}(\hat{\al}_\phi;\al)
&=E[Z_0-\hat{\be}_\phi U_0]-\al \\
&=\al+\be\th-E[\hat{\be}_\phi]\th-\al \\
&=-{\rm Bias}(\hat{\be}_\phi;\be)\th.
\end{align*}
Thus, as long as we consider the class $\hat{\al}_\phi$ as an intercept estimator, the bias reduction in intercept estimation is directly linked to that in slope estimation.
More precisely, if $\hat{\be}_\phi$ satisfies $|{\rm Bias}(\hat{\be}_\phi;\be)|\leq |{\rm Bias}(\hat{\be}^{LS};\be)|$, then $\hat{\al}_\phi$ reduces the bias of $\hat{\al}^{LS}$.

Furthermore, it is observed that
\begin{align*}
{\rm MSE}(\hat{\al}_\phi;\al)&=E[\{Z_0-\al-\be\th-(\hat{\be}_\phi U_0-\be\th)\}^2] \\
&=\tau^2+E[(\hat{\be}_\phi U_0-\be\th)^2] \\
&=\tau^2+E[\{\hat{\be}_\phi (U_0-\th)-(\hat{\be}_\phi -\be)\th\}^2] \\
&=\tau^2+E[\hat{\be}_\phi^2]\si^2+{\rm MSE}(\hat{\be}_\phi ;\be)\th^2,
\end{align*}
which implies that $\hat{\al}_\phi$ has a smaller MSE than $\hat{\al}^{LS}$ if $E[\hat{\be}_\phi^2]\leq E[(\hat{\be}^{LS})^2]$ and ${\rm MSE}(\hat{\be}_\phi;\be)\leq {\rm MSE}(\hat{\be}^{LS};\be)$.
Hence, alternative intercept estimators to $\hat{\al}^{LS}$ can be constructed from several MSE-reduced slope estimators obtained in Section \ref{sec:MSE}.

\item
If there is prior information that the slope $\be$ of \eqref{eqn:model2} lies near zero, we should positively use the prior information.
In fact, using the prior information yields a good estimator such as an admissible estimator.
See Appendix \ref{sec:B}, which discusses admissible estimation of the slope $\be$ and the intercept $\al$ under the MSE criterion.

\end{list}

\appendix
\section*{Appendix}
\section{A known variance case}\label{sec:A}

In this section, we deal with a simple case where an error variance in independent variables is known.
Here, only slope estimation is considered in model \eqref{eqn:model4}.
Assume additionally that $\bxi\ne 0_p$.

Denote the LS estimator of the slope $\be$ by $\hat\be^{LS}=\U^t\Z/\Vert\U\Vert^2$.
For the known variance case, the bias-reduced estimator \eqref{eqn:BRi} is replaced with
$$
\hat\be_\ell^{BR}=\bigg\{1+\sum_{j=1}^\ell \frac{a_j}{\Vert\U\Vert^{2j}}\bigg\}\hat{\be}^{LS},
$$
where the $a_j$ are given in \eqref{eqn:BRi} and $\ell$ is a natural number.
The following identities are needed in order to evaluate the first and second moments of $\hat\be^{LS}$ and $\hat\be_\ell^{BR}$:
\begin{align}
E\bigg[\phi(\Vert\U\Vert^2)\frac{\U^t\bxi}{\Vert\U\Vert^2}\bigg]
&=\sum_{k=0}^\infi \frac{2\la}{p+2k}P_\la(k)\int_0^\infi \phi(w)g_{p+2k}(w)\dd w, 
\label{eqn:A2}\\
E\bigg[\phi(\Vert\U\Vert^2)\frac{(\U^t\bxi)^2}{\Vert\U\Vert^{4}}\bigg]
&=\sum_{k=0}^\infi \frac{2\la(1+2k)}{p+2k}P_\la(k)\int_0^\infi \frac{\phi(w)}{w}g_{p+2k}(w)\dd w,
\label{eqn:A3}
\end{align}
where $\phi$ is a function on the positive real line and $P_\la(k)$ are the Poisson probabilities with mean $\la=\Vert\bxi\Vert^2/2$.
Identities \eqref{eqn:A2} and \eqref{eqn:A3} can be shown by using the same arguments as in the proof of Lemma \ref{lem:expectations}.

A straightforward application of identity \eqref{eqn:A2} with Lemma \ref{lem:chi} gives that
\begin{align*}
{\rm Bias}(\hat\be_\ell^{BR};\be)&=E\bigg[\bigg\{1+\sum_{j=1}^\ell \frac{a_j}{\Vert\U\Vert^{2j}}\bigg\}\frac{\U^t\bxi}{\Vert\U\Vert^2}\bigg]\be-\be \\
&=E\bigg[\frac{2\la}{p+2K}+\frac{2\la}{p+2K}\sum_{j=1}^\ell \prod_{i=1}^{j}\frac{p-2i}{p+2K-2i}\bigg]\be-\be,
\end{align*}
where $K$ is the Poisson random variable with mean $\la=\Vert\bxi\Vert^2/2$.
The same arguments as in the proof of Lemma \ref{lem:bias-BRi} yields that for $\ell\geq 1$
\begin{align*}
{\rm Bias}(\hat\be_\ell^{BR};\be)
&=E\bigg[\frac{2K}{p+2K-2}+\frac{2K}{p+2K-2}\sum_{j=1}^\ell \prod_{i=1}^{j}\frac{p-2i}{p+2K-2i-2}\bigg]\be -\be \\
&=-E\bigg[\prod_{j=1}^{\ell+1}\frac{p-2j}{p+2K-2j}\bigg]\be.
\end{align*}
In a similar fashion, we obtain
$$
{\rm Bias}(\hat\be^{LS};\be)
=-E\bigg[\frac{p-2}{p+2K-2}\bigg]\be.
$$
Thus, it is seen that $|{\rm Bias}(\hat\be_\ell^{BR};\be)|\leq |{\rm Bias}(\hat\be^{LS};\be)|$ for $1\leq \ell<(p-2)/2$.
A general result like Theorem \ref{thm:BR} can also be derived, but is omitted.

We next consider the problem of reducing the MSE of $\hat\be^{LS}$ and $\hat\be_\ell^{BR}$ in the known variance case.
Define a class of estimators as $\hat\be_\psi=\psi(\Vert\U\Vert^2)\hat\be^{LS}$, where $\psi$ is a function on the positive real line.
Taking expectation with respect to $\Z\sim\Nc_p(\be\bxi,\tau^2I_p)$ for ${\rm MSE}(\hat{\be}_\psi;\be)$, we can express ${\rm MSE}(\hat{\be}_\psi;\be)$ as
\begin{equation}\label{eqn:A4}
{\rm MSE}(\hat{\be}_\psi;\be)
=\tau^2E\bigg[\frac{\psi^2(\Vert\U\Vert^2)}{\Vert\U\Vert^2}\bigg]
 +\be^2E\bigg[\bigg\{\psi(\Vert\U\Vert^2)\frac{\U^t\bxi}{\Vert\U\Vert^2}-1\bigg\}^2\bigg].
\end{equation}

Let $W=\Vert\U\Vert^2$.
Consider a slope estimator of the form $\hat{\be}_{\bar{\psi}}=\bar{\psi}(W)\hat{\be}^{LS}$, where $\bar{\psi}(w)$ is a function of $w$ on the positive real line.
Assume that the second moment of $\hat{\be}_{\bar{\psi}}$ is finite.
Using \eqref{eqn:A2}, \eqref{eqn:A3} and \eqref{eqn:A4} leads to
\begin{align*}
&{\rm MSE}(\hat{\be}_{\psi};\be)-{\rm MSE}(\hat{\be}_{\bar{\psi}};\be) \\
&=\tau^2E\bigg[\frac{\{\psi(W)\}^2-\{\bar{\psi}(W)\}^2}{W}\bigg] \\
&\qquad +\be^2\sum_{k=0}^\infi P_\la(k)\frac{2\la}{p+2k}\int_0^\infi \De_k(w|\psi,\bar{\psi})\frac{f_{p+2k}(w)}{w}\dd w,
\end{align*}
where $
\De_k(w|\psi,\bar{\psi})=(1+2k)[\{\psi(w)\}^2-\{\bar{\psi}(w)\}^2]-2w\{\psi(w)-\bar{\psi}(w)\}$.
If $\{\psi(w)\}^2\leq \{\bar{\psi}(w)\}^2$, then $\De_k(w|\psi,\bar{\psi})\leq \De_0(w|\psi,\bar{\psi})$.
Hence, if $\{\psi(w)\}^2\leq \{\bar{\psi}(w)\}^2$ and $\De_0(w|\psi,\bar{\psi})=\{\psi(w)\}^2-\{\bar{\psi}(w)\}^2-2w\{\psi(w)-\bar{\psi}(w)\}\leq 0$, the MSE of $\hat\be_\psi$ is smaller than that of $\hat\be_{\bar\psi}$.

For a simple example, let us define $\psi_0(w)=\max[\,0,\ \min\{\bar{\psi}(w),\ 2w-\bar{\psi}(w)\}]$.
Then the resulting estimator $\hat\be_{\psi_0}=\psi_0(W)\hat\be^{LS}$ has a smaller MSE than $\hat\be_{\bar\psi}$.
By virtue of this result, we can improve on the MSEs of $\hat\be^{LS}$ and $\hat\be_\ell^{BR}$, but the details are omitted.

\section{Admissible estimators}\label{sec:B}

In this section, we present an admissible estimator of the slope $\be$ associated with proper prior distributions.
To this end, the MSE criterion is used, which means that a loss function is squared loss
\begin{equation}\label{eqn:loss}
L(\hat{\be},\be)=(\hat{\be}-\be)^2,
\end{equation}
where $\hat{\be}$ is an estimator of $\be$.
Moreover, an admissible estimator of the intercept $\al$ is derived on the basis of the admissible estimator of $\be$.

\subsection{Slope estimation}

Let $\eta=\si^{-2}$ and $\ka=\si^2/\tau^2$.
Suppose that prior densities of $\al$, $\be$, $\th$, $\bxi$ and $\eta$ are, respectively,
\begin{equation}\label{eqn:pr_al}
\pi(\al|\be,\eta,\ka)=({\rm const.})\times\Big(\frac{\ka\eta}{h_{\ka,\be}}\Big)^{1/2} \exp\Big(-\frac{c_1\ka\eta}{2h_{\ka,\be}}\al^2\Big),\quad -\infi<\al<\infi, 
\end{equation}
\begin{equation}\label{eqn:pr_be}
\pi(\be|\ka) =({\rm const.})\times\ka^{1/2} h_{\ka,\be}^{-(p+m+1)/2},\quad -\infi<\be<\infi, 
\end{equation}
\begin{equation}\label{eqn:pr_th}
\pi(\th|\eta)=({\rm const.})\times\eta^{1/2} \exp\Big(-\frac{c_2\eta}{2}\th^2\Big),\quad -\infi<\th<\infi, 
\end{equation}
\begin{equation}\label{eqn:pr_xi}
\pi(\bxi|\eta)=({\rm const.})\times\eta^{p/2} \exp\Big(-\frac{c_2\eta}{2}\Vert\bxi\Vert^2\Big),\quad \bxi\in\Rb^p, 
\end{equation}
\begin{equation}\label{eqn:pr_eta}
\pi(\eta|\be,\ka)=({\rm const.})\times h_{\ka,\be}^{-c_3/2} \eta^{c_3/2-1}\exp\Big(-\frac{\eta}{2h_{\ka,\be}}\Big),\quad \eta>0, 
\end{equation}
where $c_1$, $c_2$ and $c_3$ are certain positive constants and $h_{\ka,\be}=1+c_2+\ka\be^2$.
Suppose also that $\pi_0(\ka)$ is a suitable prior density of $\ka$ on the positive real line.
The joint prior density of $(\al,\be,\th,\bxi,\eta,\ka)$ is then proportional to
$$
\pi(\al,\be,\th,\bxi,\eta,\ka)\propto \pi_0(\ka)\ka\eta^{(p+c_3)/2}h_{\ka,\be}^{-(p+m+2+c_3)/2} e^{-(\eta/2)G_\pi},
$$
where $G_\pi=c_1\ka\al^2/h_{\ka,\be}+c_2\th^2+c_2\Vert\bxi\Vert^2+1/h_{\ka,\be}$.

The Bayes estimator of $\be$ with respect to loss \eqref{eqn:loss} is equal to a posterior mean, which has the form
$$
\hat{\be}^{PB}=\frac{\int\be\pi(\al,\be,\th,\bxi,\eta,\ka|D)\dd\al\dd\be\dd\th\dd\bxi\dd\eta\dd\ka}{\int\pi(\al,\be,\th,\bxi,\eta,\ka|D)\dd\al\dd\be\dd\th\dd\bxi\dd\eta\dd\ka},
$$
where $\pi(\al,\be,\th,\bxi,\eta,\ka|D)$ is a posterior density of $(\al,\be,\th,\bxi,\eta,\ka)$ given $D=(Z_0,\Z,U_0,\U,S)$.

\begin{lem}
If $\int_0^\infi \pi_0(\ka)\dd\ka<\infi$, then $\hat{\be}^{PB}$ can be expressed explicitly as
\begin{equation}\label{eqn:PB}
\hat{\be}^{PB}=\frac{\U^t\Z+d_1U_0Z_0}{\Vert\U\Vert^2+S+d_2U_0^2},
\end{equation}
where $d_1=c_1/(1+c_1+c_2)$ and $d_2=(c_1+c_2)/(1+c_1+c_2)$.
\end{lem}

{\bf Proof.}\ \ 
The likelihood of $(Z_0,\Z,U_0,\U,S)=(z_0,\z,u_0,\u,s)$ is written as
$$
L(z_0,\z,u_0,\u,s|\al,\be,\th,\bxi,\eta,\ka)={\rm(const.)}\times \ka^{(p+1)/2}\eta^{(2p+m+2)/2}s^{m/2-1}e^{-(\eta/2)G_L},
$$
where $G_L=\ka(z_0-\al-\be\th)^2+\ka\Vert\z-\be\bxi\Vert^2+(u_0-\th)^2+\Vert\u-\bxi\Vert^2+s$, so that the joint posterior density of $(\al,\be,\th,\bxi,\eta,\ka)$ given the data $D=(z_0,\z,u_0,\u,s)$ is expressed by
\begin{align}\label{eqn:posterior}
\pi(\al,\be,\th,\bxi,\eta,\ka|D)
&=L(z_0,\z,u_0,\u,s|\al,\be,\th,\bxi,\eta,\ka)\times\pi(\al,\be,\th,\bxi,\eta,\ka)\non\\
&\propto \pi_0(\ka)\ka^{(p+3)/2}\eta^{(3p+m+2+c_3)/2}h_{\ka,\be}^{-(p+m+2+c_3)/2}e^{-(\eta/2)G},
\end{align}
where $G=G_\pi+G_L$.
For $G$, we complete the squares with respect to $\th$, $\bxi$, $\al$ and $\be$, and then
\begin{align}\label{eqn:G}
G&=h_{\ka,\be}[\th-\{\ka\be (z_0-\al)+u_0\}/h_{\ka,\be}]^2+h_{\ka,\be}\Vert\bxi-(\ka\be\z+\u)/h_{\ka,\be}\Vert^2 \non\\
&\quad+\frac{\ka (1+c_1+c_2)}{h_{\ka,\be}}\Big\{\al-\frac{(1+c_2)z_0-\be u_0}{1+c_1+c_2}\Big\}^2 \non\\
&\quad +\frac{\ka g_0}{h_{\ka,\be}}(\be-\hat{\be}^{PB})^2 
 +\frac{\ka g_1+g_2+1}{h_{\ka,\be}},
\end{align}
where
\begin{align*}
&\hat{\be}^{PB}=\frac{\u^t\z+d_1u_0z_0}{g_0},\qquad g_0=\Vert\u\Vert^2+s+d_2u_0^2,\\
&g_1=(1+c_2)\Vert\z\Vert^2+\frac{c_1(1+c_2)z_0^2}{1+c_1+c_2}-(\hat{\be}^{PB})^2g_0, \\
&g_2=c_2u_0^2+c_2\Vert\u\Vert^2+(1+c_2)s.
\end{align*}
It thus follows that
\begin{align*}
\pi(\be,\ka|D)&\propto\int\pi(\al,\be,\th,\bxi,\eta,\ka|D)\dd\al\dd\th\dd\bxi\dd\eta\\
&\propto\pi_0(\ka)\ka^{(p+2)/2}\big\{\ka (\be-\hat{\be}^{PB})^2g_0 
 +\ka g_1+g_2+1\big\}^{-(2p+m+2+c_3)/2}.
\end{align*}
Since $\pi(\be,\ka|D)$ is symmetric at $\be=\hat{\be}^{PB}$, the Bayes estimator of $\be$ is equal to $\hat{\be}^{PB}$.

The posterior density of $\ka$ becomes
\begin{align*}
\pi(\ka|D)&\propto\int_{-\infi}^\infi\pi(\be,\ka|D)\dd\be \\
&\propto\pi_0(\ka)\ka^{(p+1)/2}(\ka g_1+g_2+1)^{-(2p+m+1+c_3)/2}.
\end{align*}
It turns out that
$$
\int_0^\infi\pi(\ka|D)\dd\ka<g_1^{-(p+1)/2}\int_0^\infi\pi_0(\ka)\dd\ka,
$$
which implies that the finiteness of $\pi(\ka|D)$ follows if $\int_0^\infi\pi_0(\ka)\dd\ka<\infi$.
Hence the proof is complete.
\qed

\begin{thm}
Assume that $\int_0^\infi\pi_0(\ka)\dd\ka<\infi$.
If $\int_0^\infi \ka^{-1}\pi_0(\ka)\dd\ka<\infi$, then $\hat{\be}^{PB}$ is admissible relative to loss \eqref{eqn:loss}.
\end{thm}

{\bf Proof.}\ \ 
When $\int_0^\infi\pi_0(\ka)\dd\ka<\infi$, $\hat{\be}^{PB}$ is proper Bayes.
Hence the admissiblity of $\hat{\be}^{PB}$ follows if the Bayes risk in terms of $\hat{\be}^{PB}$ is finite, namely
$$
\int {\rm MSE}(\hat{\be}^{PB};\be)\pi(\al,\be,\th,\bxi,\eta,\ka)\dd\al\dd\be\dd\th\dd\bxi\dd\eta\dd\ka<\infi.
$$
To prove the theorem, we shall derive a condition of the finiteness.

For real numbers $a$ and $b$ and for positive numbers $c$ and $d$, it holds true that $(a+b)^2/(c+d) \leq a^2/c+b^2/d$.
The risk of $\hat{\be}^{PB}$, namely the MSE of $\hat{\be}^{PB}$, is bounded above as
\begin{equation}\label{eqn:MSEbePB}
{\rm MSE}(\hat{\be}^{PB};\be)=E[(\hat{\be}^{PB}-\be)^2]\leq 2E[(\hat{\be}^{PB})^2]+2\be^2.
\end{equation}
Also, it is seen that
\begin{equation}\label{eqn:hatbe2-0}
\frac{(\u^t\z+d_1u_0z_0)^2}{\Vert\u\Vert^2+s+d_2u_0^2}
\leq \frac{(\u^t\z)^2}{\Vert\u\Vert^2}+\frac{(d_1u_0z_0)^2}{s+d_2u_0^2} 
\leq \Vert\z\Vert^2+\frac{d_1^2}{d_2}z_0^2
\leq \Vert\z\Vert^2+z_0^2,
\end{equation}
which implies that
\begin{equation}\label{eqn:hatbe2-1}
E[(\hat{\be}^{PB})^2]\leq E\bigg[\frac{\Vert\Z\Vert^2+Z_0^2}{\Vert\U\Vert^2+S+d_2U_0^2}\bigg].
\end{equation}
Recall that $Z_0$, $\Z$, $U_0$, $\U$ and $S$ are mutually independent.
Since $\eta(\Vert\U\Vert^2+S+U_0^2)|\th,\bxi,\eta \sim\chi_{p+m+1}^2(\de)$ with $\de=\eta(\th^2+\Vert\bxi\Vert^2)$, we observe that
\begin{align}\label{eqn:hatbe2-1d}
E\bigg[\frac{1}{\Vert\U\Vert^2+S+d_2U_0^2}\bigg]
&\leq \frac{1}{d_2}E\bigg[\frac{1}{\Vert\U\Vert^2+S+U_0^2}\bigg] \non\\
&=\frac{1}{d_2}E\bigg[\frac{\eta}{p+m-1+2K}\bigg] \leq\frac{\eta}{d_2},
\end{align}
where $K$ is the Poisson variable with mean $\eta(\th^2+\Vert\bxi\Vert^2)/2$.
It also follows that
\begin{equation}\label{eqn:hatbe2-1n}
E[\Vert\Z\Vert^2+Z_0^2]=\be^2\Vert\bxi\Vert^2+(p+1)/(\ka\eta)+(\al+\be\th)^2.
\end{equation}
Combining \eqref{eqn:hatbe2-1}, \eqref{eqn:hatbe2-1d} and \eqref{eqn:hatbe2-1n} gives that
\begin{align}\label{eqn:hatbe2-2}
E[(\hat{\be}^{PB})^2]
&\leq \frac{1}{d_2}\{\eta\be^2\Vert\bxi\Vert^2+\eta(\al+\be\th)^2+(p+1)/\ka\} \non\\
&\equiv C_1\eta\be^2\Vert\bxi\Vert^2+C_2/\ka+C_3\eta\al^2+C_4\eta\be^2\th^2+C_5\eta\al\be\th,
\end{align}
where $C_1,\ldots,C_5$ are positive constants.
Integrating both sides of \eqref{eqn:hatbe2-2} with respect to the prior densities of $\al$, $\th$ and $\bxi$, we obtain
\begin{align*}
&\int E[(\hat{\be}^{PB})^2]\pi(\al|\be,\eta,\ka)\pi(\th|\ka)\pi(\bxi|\ka)\dd\al\dd\th\dd\bxi \non\\
&\leq \frac{C_1p}{c_2}\be^2+\frac{C_2}{\ka}+\frac{C_3}{c_1}\frac{h_{\ka,\be}}{\ka}+\frac{C_4}{c_2}\be^2 \non\\
&\leq \frac{C_1p}{c_2}\frac{h_{\ka,\be}}{\ka}+\frac{C_2}{\ka}h_{\ka,\be}+\frac{C_3}{c_1}\frac{h_{\ka,\be}}{\ka}+\frac{C_4}{c_2}\frac{h_{\ka,\be}}{\ka}
\equiv \frac{C_6}{\ka} h_{\ka,\be},
\end{align*}
where $C_6$ is a positive constant.
Moreover, it follows that for a positive constant $C_7$
\begin{equation}\label{eqn:hatbe2}
\int E[(\hat{\be}^{PB})^2]\pi(\al,\be,\th,\bxi,\eta,\ka)\dd\al\dd\be\dd\th\dd\bxi\dd\eta\dd\ka
\leq C_7 \int_0^\infi \ka^{-1}\pi_0(\ka)\dd\ka
\end{equation}
because
$$
\int_{-\infi}^\infi h_{\ka,\be} \pi(\be|\ka)\dd\be=\frac{p+m-1}{p+m-2}(1+c_2).
$$

In the same way as above, taking expectation of $\be^2$ with respect to the prior densities yields that, for a positive constant $C_8$,
\begin{equation}\label{eqn:be2}
\int \be^2\pi(\al,\be,\th,\bxi,\eta,\ka)\dd\al\dd\be\dd\th\dd\bxi\dd\eta\dd\ka=C_8\int_0^\infi \ka^{-1}\pi_0(\ka)\dd\ka.
\end{equation}
By combining \eqref{eqn:MSEbePB}, \eqref{eqn:hatbe2} and \eqref{eqn:be2}, the Bayes risk of $\hat{\be}^{PB}$ can be bounded above as
$$
\int {\rm MSE}(\hat{\be}^{PB};\be)\pi(\al,\be,\th,\bxi,\eta,\ka)\dd\al\dd\be\dd\th\dd\bxi\dd\eta\dd\ka
\leq C_9\int_0^\infi \ka^{-1}\pi_0(\ka)\dd\ka 
$$
for a positive constant $C_9$.
Hence, if the r.h.s. of the above inequality is finite, the Bayes risk of $\hat{\be}^{PB}$ is finite.
\qed

\subsection{Intercept estimation}

Next, we address admissible estimation of the intercept $\al$ under the squared loss $(\hat{\al}-\al)^2$.

An admissible estimator of $\al$ is derived with the aid of proper priors \eqref{eqn:pr_al}--\eqref{eqn:pr_eta}.
Let $\pi_0(\ka)$ be a prior density of $\ka$ such that $\int_0^\infi\pi_0(\ka)\dd\ka<\infi$.
From \eqref{eqn:posterior} and \eqref{eqn:G}, we obtain the Bayes estimator, namely a posterior mean,
$$
\hat{\al}^{PM}=d_1^*Z_0-d_2^*\hat{\be}^{PM} U_0,
$$
where $d_1^*=(1+c_2)/(1+c_1+c_2)$, $d_2^*=1/(1+c_1+c_2)$ and $\hat{\be}^{PM}$ is given in \eqref{eqn:PB}.

The admissibility of $\hat{\al}^{PM}$ is based on the following theorem.

\begin{thm}
If $\int_0^\infi\ka^{-1}\pi_0(\ka)\dd\ka<\infi$ and $c_3>2$, then $\hat{\al}^{PB}$ is admissible relative to the squared loss.
\end{thm}

{\bf Proof.}\ \
The MSE of $\hat{\al}^{PM}$ is bounded above by
\begin{align}\label{eqn:MSEal}
{\rm MSE}(\hat{\al}^{PM};\al)&=E[(d_1^*Z_0-d_2^*\hat{\be}^{PM} U_0-\al)^2] \non\\
&\leq 3(d_1^*)^2E\big[Z_0^2\big]+3(d_2^*)^2E\big[\big(\hat{\be}^{PM}\big)^2 U_0^2\big]+3\al^2 \non\\
&\leq 3\{(\al+\be\th)^2+1/(\ka\eta)\}+3E\big[\big(\hat{\be}^{PM}\big)^2 U_0^2\big]+3\al^2.
\end{align}
Here, using the same arguments as in \eqref{eqn:hatbe2-0} and \eqref{eqn:hatbe2-1n} leads to
\begin{align}\label{eqn:PMU2}
E\big[\big(\hat{\be}^{PM}\big)^2 U_0^2\big]
&\leq E\bigg[\frac{\Vert\Z\Vert^2+Z_0^2}{\Vert\U\Vert^2+S+d_2U_0^2}U_0^2\bigg]
 \leq \frac{1}{d_2}E\big[\Vert\Z\Vert^2+Z_0^2\big] \non\\
&=\frac{1}{d_2}\big\{\be^2\Vert\bxi\Vert^2+(p+1)/(\ka\eta)+(\al+\be\th)^2\big\}.
\end{align}
Combining \eqref{eqn:MSEal} and \eqref{eqn:PMU2}, we can write the upper bound of ${\rm MSE}(\hat{\al}^{PM};\al)$ as
\begin{equation}\label{eqn:up_mse_al1}
{\rm MSE}(\hat{\al}^{PM};\al)\leq C_1^*\be^2\Vert\bxi\Vert^2+C_2^*/(\ka\eta)+C_3^*\al^2+C_4^*\be^2\th^2+C_5^*\al\be\th,
\end{equation}
where $C_1^*,\ldots,C_5^*$ are positive constants.
Taking expectation of \eqref{eqn:up_mse_al1} with respect to \eqref{eqn:pr_al}, \eqref{eqn:pr_th} and \eqref{eqn:pr_xi}, we obtain
\begin{align}\label{eqn:up_mse_al2}
&\int {\rm MSE}(\hat{\al}^{PM};\al)\pi(\al|\be,\eta,\ka)\pi(\th|\eta)\pi(\bxi|\eta)\dd\al\dd\th\dd\bxi \non\\
&=(pC_1^*+C_4^*)\be^2/(c_2\eta)+C_2^*/(\ka\eta)+C_3^*h_{\ka,\be}/(c_1\ka\eta) \non\\
&\leq (pC_1^*+C_4^*)h_{\ka,\be}/(c_2\ka\eta)+C_2^*h_{\ka,\be}/(\ka\eta)+C_3^*h_{\ka,\be}/(c_1\ka\eta) \non\\
&= (pC_1^*/c_2+C_2^*+C_3^*/c_1+C_4^*/c_2)h_{\ka,\be}/(\ka\eta).
\end{align}
Next, taking expectation of \eqref{eqn:up_mse_al2} with respect to \eqref{eqn:pr_eta} gives that for $c_3>2$
\begin{equation*}
\int {\rm MSE}(\hat{\al}^{PM};\al)\pi(\al|\be,\eta,\ka)\pi(\th|\eta)\pi(\bxi|\eta)\pi(\eta|\be,\ka)\dd\al\dd\th\dd\bxi\dd\eta
\leq C_6^*/\ka.
\end{equation*}
where $C_6^*=(pC_1^*/c_2+C_2^*+C_3^*/c_1+C_4^*/c_2)/(c_3-2)$.
Hence we obtain
$$
\int {\rm MSE}(\hat{\al}^{PM};\al)\pi(\al,\be,\th,\bxi,\eta,\ka)\dd\al\dd\be\dd\th\dd\bxi\dd\eta\dd\ka
\leq C_6^*\int_0^\infi\ka^{-1}\pi_0(\ka)\dd\ka,
$$
which complete the proof.
\qed



\end{document}